\documentclass[10pt]{article}
\usepackage{amsmath,amssymb,amsthm,morefloats}
\usepackage{epsfig,psfrag,url}
\usepackage{graphicx}
\usepackage{verbatim}
\usepackage{subfigure}

\numberwithin{equation}{section}

\newcommand{\goto}{\rightarrow}
\newcommand{\bigo}{{\mathcal O}}

\def\Xint#1{\mathchoice
   {\XXint\displaystyle\textstyle{#1}}%
   {\XXint\textstyle\scriptstyle{#1}}%
   {\XXint\scriptstyle\scriptscriptstyle{#1}}%
   {\XXint\scriptscriptstyle\scriptscriptstyle{#1}}%
   \!\int}
\def\XXint#1#2#3{{\setbox0=\hbox{$#1{#2#3}{\int}$}
     \vcenter{\hbox{$#2#3$}}\kern-.5\wd0}}

\def\dashint{\Xint-}
\DeclareMathOperator{\sech}{sech}

\DeclareMathOperator{\diag}{diag}
\DeclareMathOperator{\res}{Res}
\DeclareMathOperator{\imag}{Im}

\DeclareMathOperator{\id}{id}
\DeclareMathOperator{\sign}{sign}
\DeclareMathOperator{\trace}{tr}

\newenvironment{choices}{\left\{ \begin{array}{ll}}{\end{array}\right.}
\newcommand\when{&\text{if~}}
\newcommand\otherwise{&\text{otherwise}}

\newenvironment{mat}{\left[\begin{array}{ccccccccccccccc}}{\end{array}\right]}
\newcommand\bcm{\begin{mat}}
\newcommand\ecm{\end{mat}}

\newcommand{\bea}{\begin{eqnarray}}
\newcommand{\eea}{\end{eqnarray}}
\newcommand{\bean}{\begin{eqnarray*}}
\newcommand{\eean}{\end{eqnarray*}}
\newcommand{\ba}{\begin{array}}
\newcommand{\ea}{\end{array}}
\newcommand{\beqs}{\begin{equation*}\begin{split}}

\newtheorem{example}{Example}[section]

\newtheorem{remark}{Remark}[section]

\newtheorem{lemma}{Lemma}[section]

\newtheorem{theorem}{Theorem}[section]

\long\def\symbolfootnote[#1]#2{\begingroup%
\def\thefootnote{\fnsymbol{footnote}}\footnote[#1]{#2}\endgroup}

\begin{document}
\title{On the application of GMRES to oscillatory singular integral equations}
\author{Thomas Trogdon$^1$\\
\phantom{.}\\
Courant Institute of Mathematical Sciences\\
 New York University\\
251 Mercer St.\\
New York, NY 10012, USA \\}
\maketitle

\footnotetext[1]{Email: trogdon@cims.nyu.edu}

\begin{abstract}

We present a new method for the numerical solution of singular integral equations on the real axis.  The method's value stems from a new formula for the Cauchy integral of  a rational function with an oscillatory exponential factor. The inner product of such functions is also computed explicitly.  With these tools in hand, the GMRES algorithm is applied to both non-oscillatory and oscillatory singular integral equations.  In specific cases, ideas from Fredholm theory and Riemann--Hilbert problems are used to motivate preconditioners for these singular integral equations.  A significant acceleration in convergence is realized for these examples.  This presents a useful link between the theory of singular integral equations and the numerical analysis of such equations.  Furthermore, this method presents a first step towards a solver for the inverse scattering transform that does not require the deformation of a Riemann--Hilbert problem.
\end{abstract}

\section{Introduction}

The numerical analysis of singular integral equations has historically been an important topic \cite{mikhlin,prossdorf}.  Such equations are often difficult to analyze numerically because operators involved are not compact.  Consider the operator equation
\begin{align*}
(\id + \mathcal T)u = f,
\end{align*}
where $\id$ is the identity operator and $\mathcal T$ is a non-compact operator.  If $\mathcal T$ is replaced with any finite-dimensional approximation $\mathcal T_n$, the lack of compactness guarantees that $\|\mathcal T_n - \mathcal T\| \not \goto 0$ as $n \goto \infty$.  This fact complicates both the construction of numerical schemes and corresponding proofs of convergence \cite{prossdorf}.  The main goal of this paper is to develop an approach that avoids a finite-dimensional approximation of a non-compact oscillatory singular integral operator by computing
\begin{align*}
\frac{1}{2 \pi i} \int_{\mathbb R} \left[ \left( \frac{ x - \beta i}{x + \beta i} \right)^j -1 \right] e^{i \alpha x} \frac{dx}{x-z},
\end{align*}
explicitly for all $j \in \mathbb Z$, $\alpha \in \mathbb R$ and $\beta > 0$.  The GMRES \cite{GMRES-original} algorithm is applied once the action of the relevant operator is computed exactly.

  In the previous two decades there has been increased interest in the solution of singular integral equations largely due to their connections with Riemann--Hilbert problems (RHPs) \cite{AblowitzClarksonSolitons,BealsCoifman,DeiftOrthogonalPolynomials}. A numerical approach can be found in \cite{DienstfreyThesis} where the author concentrates on the solution of an RHP that arises in the computation of specific nonlinear special functions.  An expanded treatment was developed in \cite{SOHilbertTransform,SOPainleveII,SORHFramework} where singular integral equations on fairly general domains are solved.  Again, this method is applicable to the computation of many nonlinear special functions including solutions of partial differential equations.  See \cite{TrogdonThesis} for a demonstration of the wide array of functions that can be computed with this method.

The prototypical domains for singular integral equations are  $\mathbb R$, the unit circle $\{|z| = 1\}$ and the interval $[-1,1]$ \cite{prossdorf}.  There are many reasons to consider singular integral equations posed on more general domains.  A domain that consists of many contours that intersect at the origin arises in inverse scattering for higher-order systems \cite{Ontheline}.  Techniques in the recent method of Fokas lead to similar domains \cite{FokasUnified}.  Even when the equation is initially posed on a simple domain such as $\mathbb R$, contour deformations (in the case of RHPs) that reduce oscillation to exponential decay transform $\mathbb R$ to a more complicated domain \cite{DeiftZhouAMS}.  When this problem is approached from a numerical point of view, similar contours are found \cite{TrogdonSONNSD,TrogdonSORMT,TrogdonSONLS,TrogdonSOKdV}.  It is important to understand if these contour deformations are necessary for the numerical analysis of the problem.

In this paper, we present what we believe is the first-known method for the numerical solution of oscillatory singular integral equations on the line without any contour deformation.  Due to the domain being restricted to the real line, it is clear that our method is in not a replacement for \cite{SORHFramework}.  Despite this, the method is the first step towards a general framework for oscillatory singular integral equations that would encompass the integral equations obtained in the inverse scattering transform \cite{DeiftZhouAMS}.  Furthermore, we treat singular integral equations that have slow decay in the coefficient functions at infinity.

As noted above, the GMRES algorithm is applied to solve singular integral equations. In examples, we consider two choices of preconditioners to accelerate convergence of GMRES.  The first choice is the so-called Fredholm regulator.  For a Fredholm operator $\id + \mathcal T$, the Fredholm regulator $\mathcal R$ is an operator chosen so that $\mathcal R(\id + \mathcal T) = \id + \mathcal K$ where $\mathcal K$ is a compact operator.  Empirically, GMRES applied to $\mathcal R(\id + \mathcal T)u = \mathcal R f$ converges faster than it does when applied to $(\id + \mathcal T)u = f$.  The second choice for a preconditioner (see Section~\ref{Section:OscLine}) is motivated directly by the method of nonlinear steepest descent \cite{DeiftZhouAMS}.  A significant speedup in convergence is realized.  We see that preconditioners can be motivated by the theory of singular integral equations.

In this paper, four main conclusions are reached:
\begin{itemize}
\item Our formula for the Cauchy integral of rational functions with an oscillatory exponential factor has wide applicability,
\item GMRES is an effective tool for the numerical solution of oscillatory singular integral equations on the real axis,
\item in specific examples preconditioning operators, which accelerate the convergence of GMRES, can be motivated from the underlying singular integral equation theory, and
\item there is a strong indication that GMRES should be used in future research on oscillatory singular integral equations.
\end{itemize}

The paper is separated into two parts.  The first part (Sections~\ref{Section:non-osc}--\ref{Section:NonOscExamples}) is concerned with the development of the tools in the absence of oscillations.  The reasoning for this is two-fold.  First, this allows the demonstration of the broad ideas of the paper with out extra complication.  Second, we discuss how the method breaks down when oscillations are introduced (see Remark~\ref{non-osc-issues}) .  This motivates the further developments that follow.


In the second part of the paper, (Sections~\ref{Section:OscBasis} and \ref{Section:OscLine}) a new formulae for the Cauchy integral acting on basis of oscillatory rational functions is derived.  This allows the accurate computation of the Cauchy integral of a Fourier-type integrand.    We also present new formulae for the inner product of these oscillatory basis functions and this is used for oscillatory quadrature.  Similar results exist in the literature (see \cite{KellerPractical,WangHilbert}) for different bases and different integration domains.  Using these ideas, the Fourier transforms of functions that decay slowly can be computed.  In addition, we show an application of the method to the solution of linear partial differential equations.  Finally, all the results in Section~\ref{Section:OscLine} are combined to solve oscillatory singular integral equations stably for all parameter values, after preconditioning.

\section{Non-oscillatory basis functions}\label{Section:non-osc}

As a motivating problem, we consider the problem of computing the Cauchy integral
\begin{align}\label{Cauchy-Integral}
\mathcal C_{\Gamma} f(z) = \frac{1}{2 \pi i} \int_{\Gamma} \frac{f(s)}{s-z} ds.
\end{align}
Define $\mathbb U = \{|z| = 1\}$ with counter-clockwise orientation.  In one approach (see, for example, \cite{SOHilbertTransform}) a uniform approximation of $f(s)$ in a Laurent series is found:
\begin{align}\label{f-uniform}
\sup_{s \in \mathbb U}|f(s)-f_n(s)| \goto 0 \mbox{ as } n \goto \infty, \quad f_n(s) = \sum_{j=-n}^n a_{j,n} s^j.
\end{align}
It follows from straightforward contour integration that
\begin{align*}
\mathcal C_{\mathbb U} f(z) \approx \begin{choices} \displaystyle \sum_{i=0}^n a_{j,n} z^j, \when |z| < 1,\\
\displaystyle-\sum_{j=-n}^{-1} a_{j,n} z^j, \when |z| > 1, \end{choices}
\end{align*}
is uniformly accurate on sets bounded away from $\mathbb U$.  An efficient way of computing the approximations $f_n$ of $f$ is found through the fast Fourier transform (see Example \ref{Example:RationalApprox}).  We use this idea to motivate a method for computing $\mathcal C_{\mathbb R}$.

Consider the family of M\"obius transformations
\begin{align*}
M_\beta(z) = \frac{z - i \beta}{z+i\beta}, ~~ M^{-1}_\beta(z) = \frac{\beta}{i} \frac{z+1}{z-1},~~ \beta >0.
\end{align*}
Each of these transformations $M_\beta$ maps the real axis to the unit circle.  We look for an approximation $g_n$ of $g = f \circ M_\beta^{-1}$ in terms of a Laurent series since $g$ is defined on the unit circle.  Regularity conditions on $f$ must be imposed.  For our purposes we require that $f$ is smooth and decay rapidly at infinity.  The decay requirement can be relaxed provided $f$ has smoothness at infinity on the Riemann sphere.  Assume we have a uniform approximation of $g$:
\begin{align}\label{laurent-approx}
g(s) \approx g_n(s) = \sum_{j=-n}^n a_{j,n} s^j.
\end{align}
It follows that for $\beta > 0$ \cite{SOHilbertTransform},
\begin{align}\label{line-formula}
\mathcal C_{\mathbb R} f(z) \approx \begin{choices} \displaystyle \sum_{j=0}^n a_{j,n} M_\beta^j(z) - \sum_{j=0}^n a_{j,n}, \when \imag z > 0,\\
\displaystyle\sum_{j=-n}^{-1} a_{j,n} M_\beta^j(z) - \sum_{j=-n}^{-1} a_{j,n}, \when \imag z < 0,
\end{choices}
\end{align}
and the approximation is uniformly accurate on sets bounded away from $\mathbb R$.  Again, the coefficients $a_{j,n}$ may be approximated well with the fast Fourier transform (see Example~\ref{Example:RationalApprox}).

The range of the Cauchy integral when acting on smooth, rapidly decaying functions is conveniently represented by the basis
\begin{align}\label{non-osc-basis}
\left\{R_j(z)\right\}_{j=-\infty}^\infty, ~~ R_j(z) = M^j_\beta(z) -1, ~~\beta >0.
\end{align}
We suppress the dependence on $\beta>0$ when writing the basis and make any choices clear below.  We note that a minor modification of $M^j_\beta(z)$ produces the functions used for approximation of Cauchy integrals in \cite{SOHilbertTransform}.  Essentially, $R_j$ is found by subtracting the asymptotic (large $z$) behavior of $M^j_\beta(z)$ to make the resulting functions square integrable.

\begin{remark}
We use the parameter $\beta$ in anticipation of cases where this extra degree of freedom is useful.  For example, exactly expressing the function $z \mapsto (z + 2i)^{-1}$ in terms of $M_\beta(z)$ requires an infinite sum for $\beta \neq 2$.   In some applications, the correct choice of $\beta$ may result in a sparse approximation.
\end{remark}

\begin{example}[Rational approximation]\label{Example:RationalApprox}
Our task is to expand a function in the basis $\{R_j(z)\}_{j=-\infty}^\infty$.  We use ideas from \eqref{laurent-approx}.  Let $g = f \circ M_1^{-1}$ ($\beta = 1$).  Applying the fast Fourier transform to $g(e^{i\theta})$ sampled on a uniform grid we obtain
\begin{align*}
g(e^{i\theta}) \approx \sum_{n=-N}^N a_{n,N} e^{in\theta} \Leftrightarrow f(x) \approx \sum_{n=-N}^N a_{n,N} M_1^n(z). 
\end{align*}
If the grid $\{ 2\pi j/(2N+1) : j = 0,1,2,\ldots,2N\}$ is chosen then  $\sum_n a_{n,N} = 0$ because $f(\infty) = 0$ and $\infty$ is an interpolation point.  Therefore
\begin{align*}
f(x) \approx \sum_{n=-N}^N a_{n,N} R_n(z).
\end{align*}
For concreteness, consider approximating $f(x) = e^{-x^2}$.  See Figure~\ref{Gaussian-approximation} for numerical results.   The derivation of precise error bounds for this method was performed recently in \cite{Trogdon2014}.

\begin{figure}[tb]
\centering
\includegraphics[width=.5\linewidth]{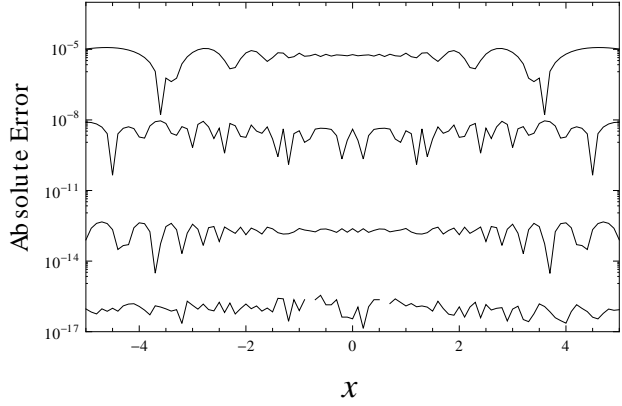}
\caption{\label{Gaussian-approximation}  The error in the approximation of $f(x) = e^{-x^2}$ in the basis $\{R_j\}_{j=-N}^N$ using the fast Fourier transform to compute the coefficients.  The figure shows the absolute error for  $N = 20,40,80$ and $160$. As expected, we see spectral convergence.}
\end{figure}

\end{example}

\subsection{Properties of the basis}

Additional practical and theoretical properties of this basis are now laid out.    Our first result anchors the theoretical developments.

\subsubsection{Density}

\begin{theorem}
$\left\{R_j(z)\right\}_{j=-\infty}^\infty$ is a basis for $L^2(\mathbb R)$.
\begin{proof}
First, it is clear that $R_j(z)$ is square integrable for all $j$. For $j > 0$, we write
\begin{align*}
R_j(z) &= \left(\frac{z-i\beta}{z+ i\beta}\right)^j - \left(\frac{z+i\beta}{z+ i\beta}\right)^j = \frac{P_{j}(z)}{(z+i\beta)^j},\\
P_{j}(z) &= \sum_{k=1}^{j} \left(\begin{array}{c} j \\ k \end{array}\right) (z+i\beta)^{j-k} (-2i\beta)^{k},
\end{align*}
from the Binomial Theorem.  More abstractly,
\begin{align*}
P_{j} (z) = \sum_{k=1}^{j} c_{k,j} (z+i\beta)^{j-k},
\end{align*}
for some ($\beta$-dependent) coefficients $c_{k,j}$ so that
\begin{align}\label{rat-sum}
R_j(z) = \sum_{k=1}^j \frac{c_{k,j}}{(z+i\beta)^k}.
\end{align}
It is well-known that the set
\begin{align*}
\left\{ K_j(z) \right\}_{j=-\infty}^\infty,~~ K_j(z) = \frac{1}{z+i \beta} \left(\frac{z-i \beta}{z+i\beta} \right)^{j-1},
\end{align*}
forms an orthogonal basis for $L^2(\mathbb R)$ \cite[p.~195]{SteinRA}. For $j > 0$, each element of this basis can be expressed in the form of \eqref{rat-sum} showing that there is a 1-1 correspondence between the two bases. For $j < 0$ we have to work a bit harder.  The above arguments show for some new coefficients $\tilde c_{j,k}$
\begin{align}\label{rat-sum-neg}
R_j(z) = \sum_{k=1}^{-j} \frac{\tilde c_{k,j}}{(z-i\beta)^k}, ~~ j < 0.
\end{align}
For $j \leq 0$
\begin{align*}
K_j(z) = \frac{z-i\beta}{z+i\beta} \sum_{k=1}^{-j} \frac{d_{k,j}}{(z-i\beta)^k},
\end{align*}
for some coefficients $d_{k,j}$.  This shows that
\begin{align*}
\frac{z+i\beta}{z-i\beta}K_j(z)
\end{align*}
can be expressed in terms of $R_j(z)$.  Define an invertible, bounded linear transformation $\mathcal T$ on $L^2(\mathbb R)$ by
\begin{align*}
\mathcal T\left(\sum_{j=-n}^n a_j K_j(z)\right) = \sum_{j=1}^n a_j K_j(z) +  \sum_{j=-n}^0 a_j \frac{z+i\beta}{z-i\beta} K_j(z), ~~ \text{ for all } n> 0.
\end{align*}
For any function $g \in L^2(\mathbb R)$ expand
\begin{align*}
\sum_{j=-n}^n a_j K_j(z) \goto \mathcal T^{-1}g \text{ in } L^2(\mathbb R) ~~\text{ so that }~~ \sum_{j=-n}^n a_j \mathcal TK_j(z) \goto g \text{ in } L^2(\mathbb R).
\end{align*}
It follows that $\sum_{j=-n}^n a_j \mathcal TK_j(z)$ can be expressed in terms of $R_j$.  This proves the theorem.
\end{proof}
\end{theorem}

\subsubsection{Action of the Cauchy operators}

As discussed in the introduction, the main motivation for considering the Cauchy integral is to compute the Cauchy operators.  For our purposes here, the Cauchy operators are defined by
\begin{align*}
\mathcal C_{\mathbb R}^\pm f(x) = \lim_{\epsilon \goto 0^+} \frac{1}{2 \pi i} \int_{\mathbb R} \frac{f(s)}{s-(x\pm i \epsilon)} ds.
\end{align*}
It is known that if $f \in L^2(\mathbb R)$ then this limit exists a.e. and is an $L^2(\mathbb R)$ function that satisfies $\| \mathcal C_{\mathbb R}^\pm f\|_{L^2(\mathbb R)} \leq \| f\|_{L^2(\mathbb R)}$ \cite{DeiftOrthogonalPolynomials}.  The so-called Plemelj Lemma also holds:
\begin{align}\label{Plemelj}
\mathcal C_{\mathbb R}^+ f - \mathcal C_{\mathbb R}^- f = f. 
\end{align}
Straightforward contour integration shows:
\begin{itemize}
\item $\mathcal C^+_{\mathbb R} R_j(z) = R_j(z)$ for $j >0$,
\item $\mathcal C^-_{\mathbb R} R_j(z) = 0$ for $j >0$,
\item $\mathcal C^-_{\mathbb R} R_j(z) = -R_j(z)$ for $j <0$, and
\item $\mathcal C^+_{\mathbb R} R_j(z) = 0$ for $j <0$.
\end{itemize}

\subsection{Multiplication}
We address a practical question concerning the multiplication of two functions expanded in the basis $\{R_j(z)\}$.  A useful identity is
\begin{align}
R_j(z)R_k(z) &= (M_\beta^j(z) -1)(M_\beta^k(z)-1) = M_\beta^{j+k}(z) - M_\beta^k(z) -M_\beta^j(z) + 1\notag\\
&= (M_\beta^{j+k}(z)-1) - (M_\beta^j(z) - 1) - (M_\beta^j(z) - 1) = R_{k+j}(z) -R_j(z)-R_k(z).\label{non-osc-mult}
\end{align}
Multiplication acts almost as it does on exponentials.  Using this identity we consider for $m > n$
\begin{align*}
\left(\sum_{j=-m}^m a_j R_j(z) \right) \left(\sum_{k=-n}^n b_k R_k(z) \right) &= \sum_{-n \leq j,k \leq n} a_jb_k ( R_{k+j}(z) -R_j(z)-R_k(z))\\
&+ \sum_{j \geq m} \sum_{k =1}^n a_jb_k ( R_{k+j}(z) -R_j(z)-R_k(z))\\
&+ \sum_{j \leq m} \sum_{k =1}^n a_jb_k ( R_{k+j}(z) -R_j(z)-R_k(z))\\
&=  \sum_{l=-m+n}^{m-n} \left( \sum_{k=-n}^n a_jb_{l-j} \right) R_l(z) + \sum_{l=-m-n}^{-m+n-l} \left( \sum_{k=l+m}^n a_jb_{l-j} \right) R_l(z)\\
&+ \sum_{l=m-n+1}^{m+n} \left( \sum_{k=-n}^{l-m} a_jb_{l-j} \right) R_l(z) + \sum_{l=-m}^m a_l \left( \sum_{k=-n}^n b_k \right) R_l(z) \\
&+ \sum_{l=-n}^n \left( \sum_{j=-m}^m a_j \right) b_l R_l(z).
\end{align*}

\begin{remark}
We also consider the approximation of matrix-valued functions.  In this case the coefficients $a_j$ and $b_k$ are matrices and the same multiplication formula holds when it is taken into account that $a_j$ and $b_k$ do not necessarily commute.
\end{remark} 

\section{An integration formula}\label{Appendix:Quadrature}

Provided that $f \in L^1(\mathbb R)$ it follows that \cite{TrogdonSOKdV}
\begin{align*}
-2\pi i\lim_{z \goto \infty} z \mathcal C_{\mathbb R} f(z) = \int_{\mathbb R} f(x) dx.
\end{align*}
It is easy to see that
\begin{align}\label{res-inf}
\lim_{z \goto \infty} z R_j(z) = -2 i j \beta
\end{align}
so that if
\begin{align*}
f(x) = \sum_{j = -\infty}^\infty a_j R_j(x),~ \text{ with } \sum_{j = -\infty}^\infty |j||a_j| < \infty,
\end{align*}
then
\begin{align*}
\int_{\mathbb R}f(x)dx = -4 \pi \beta \sum_{j=1}^\infty j a_j= 4 \pi \beta \sum_{j=-\infty}^{-1} j a_j = -2 \pi \beta \sum_{j\neq 0} |j| a_j.
\end{align*}
This formula is used to integrate solutions of integral equations below.

\begin{remark}
This integration formula can be interpreted as a classical quadrature rule in the following way.  If coefficients $a_j$ are approximated via the fast Fourier transform as in Example~\ref{Example:RationalApprox} then we have the formula (using $2n+1$ quadrature nodes for simplicity)
\begin{align*}
a_j &= \frac{1}{2n+1} \sum_{l=0}^{2n} e^{-ijt_\ell} f(M^{-1}_\beta(e^{i t_\ell})),\\
t_\ell &= 2 \pi \frac{\ell}{2n +1}.
\end{align*}
Therefore, with the convention that $a_0 = 0$,
\begin{align*}
\int_{\mathbb R} f(x) dx &\approx -2 \pi \beta \sum_{j=-n}^n \sum_{l=0}^{2n}\frac{|j|}{2n+1} e^{-ijt_\ell} f(M^{-1}_\beta(e^{i t_\ell})) \\
&= -\sum_{\ell=0}^{2n} f(M_\beta^{-1}(e^{it_\ell})) \left(\sum_{j = -n}^{n} \frac{2 \pi \beta |j|}{2n +1}  e^{-ijt_\ell}    \right).
\end{align*}
We obtain the quadrature nodes $\{M_\beta^{-1}(e^{i t_\ell})\}$ and weights
\begin{align*}
\omega_\ell = -\sum_{j = -n}^{n} \frac{2 \pi \beta |j|}{2n +1}  e^{-ijt_\ell} = -\frac{2 \pi \beta}{2n +1} \sum_{j=1}^n j (e^{-ijt_\ell}+ e^{ijt_\ell}) = -\frac{4 \pi \beta}{2n +1} \sum_{j=1}^n j \cos(j t_\ell).
\end{align*}
From the nature of the method one would expect that this approximation convergences spectrally fast to the integral of $f$ provided that $f$ is smooth and rapidly decaying.  Despite this, it is not clear that this method has any advantage over, say, the trapezoidal rule after a change of variables. A detailed examination of these ideas is not performed here but see \cite{Trogdon2014} for a detailed error analysis.
\end{remark}

\subsubsection{The inner product}

Another important aspect for the basis is the computation of inner products. We derive a formula for
\begin{align*}
L_{j,k} = \int_{\mathbb R} R_j(z) \overline{R_k(z)} dz
\end{align*}
with contour integration.  First note that $\overline{R_k(z)} = R_{-k}(z)$ so that we must compute
\begin{align*}
\dashint_{\mathbb R} (R_{j-k}(z) - R_j(z) - R_{-k}(z))dz.
\end{align*}
It suffices to compute the principal value integral of $R_j(z)$ for all $j$.  If $j > 0$, consider
\begin{align}\label{c-r}
\int_{-r}^r R_j(z) dz = \int_{C^+_r} R_j(z) dz,
\end{align}
where $C^+_r =\{re^{i \theta}: 0 \leq \theta \leq \pi\}$ with clockwise orientation.  As $r \goto \infty$ this converges to $-i \pi$ times the residue of $R_j(z)$ at infinity (see  \eqref{res-inf}).  For $j <0$, replace $C^+_r$ with $C^-_r = \{re^{i \theta}: -\pi \leq \theta \leq 0\}$ with counter-clockwise orientation.  Therefore
\begin{align*}
\dashint_{\mathbb R} R_j(z) dz = - 2 \pi |j| \beta.
\end{align*}
We have the general formula
\begin{align}\label{non-osc-inner}
L_{j,k} = -2 \pi \beta(|j-k| - |j| - |k|).
\end{align}
Note that if $j$ and $-k$ have the same sign then this formula implies $L_{j,k} = 0$.

\begin{remark}
The inner product needs to be generalized to matrix-valued functions.  If $f$ and $g$ are $n\times n$ matrix-valued functions the appropriate inner product is
\begin{align}\label{matrix-inner}
\langle f, g\rangle = \int_{\mathbb R}  \trace f(x)g^*(x) dx,
\end{align}
where $^*$ represents the Hermitian conjugate of the matrix.  When $f$ and $g$ are each expressed as a series in $\{R_j\}$ this inner product is computed with \eqref{non-osc-inner} in a straightforward way using the linearity of the trace operation.
\end{remark}

\section{Singular integral equations on the line}\label{Section:SIE-line}

We consider singular integral equations that arise in the solution of RHPs on the line.  In short, an RHP consists of finding a sectionally analytic function $\Phi(z)$ that satisfies
\begin{align}\label{jump-condition}
\Phi^+(x) &= \Phi^-(x) G(x) + F(x), ~~ x \in \mathbb R,\\
\Phi^\pm (x) &\triangleq \lim_{\epsilon \goto 0^+} \Phi(x\pm i\epsilon)\notag.
\end{align}
Here $G$ and $F$ are definite functions on $\mathbb R$. The function $G$ is referred to as the jump matrix.  In general, $G$ and $F$ may be $n\times n$ matrix-valued functions which forces $\Phi$ to be matrix-valued with the same dimension.  We use the normalization condition
\begin{align*}
\lim_{z \goto \infty} \Phi(z) = I,
\end{align*}
where $I$ is the $n\times n$ identity matrix.  If we impose that \eqref{jump-condition} should hold a.e. and $\Phi^\pm-I \in L^2(\mathbb R)$ then 
\begin{align}\label{ansatz}
\Phi(z) = I + \mathcal C_{\mathbb R} u(z),
\end{align}
for $u \in L^2(\mathbb R)$ \cite{TrogdonThesis}. For this purpose, $L^2(\mathbb R)$ is appropriately generalized for matrix-valued functions using the inner product \eqref{matrix-inner}.

We use this representation  to convert an RHP to a singular integral equation.  The substitution of \eqref{ansatz} along with \eqref{Plemelj} produces
\begin{align}\label{sie}
u- \mathcal C_{\mathbb R}^- u \cdot (G-I) = G-I + F.
\end{align}
In what follows we assume $G-I, F \in L^2\cap L^\infty(\mathbb R)$ and the notation $\mathcal C[G;\mathbb R]$ is used to denote the operator
\begin{align*}
u \mapsto u - \mathcal C_{\mathbb R} u \cdot (G-I).
\end{align*}
Unless otherwise noted, for an operator $\mathcal M$, $\|\mathcal M\|$ denotes the standard operator norm on $L^2(\mathbb R)$.

In our examples, we can approximate each component of $G-I$ and $F$ accurately in the basis $\{R_j\}$.  The following result justifies replacing $G-I$ and $F$ with these approximations even though our approximations are not integrable.  Only the convergence of principal-value integrals is needed.  We use $\hat G$ and $\hat F$ to denote approximations of $G$ and $F$, respectively.
\begin{lemma}
Assume that $\mathcal C[G;\mathbb R]$ is invertible on $L^2(\mathbb R)$.  For $0 < \epsilon < 1/\|\mathcal C[G;\mathbb R]^{-1}\|$ assume  $\|\hat G - G\|_{L^2\cap L^\infty(\mathbb R)} < \epsilon$ and $\|\hat F - F\|_{L^2\cap L^\infty(\mathbb R)} < \epsilon$  then $\mathcal C[\hat G;\mathbb R]$ is also invertible and
\begin{align*}
\|\mathcal C[G;\mathbb R]^{-1} - \mathcal C[\hat G;\mathbb R]^{-1}\|\leq  \epsilon \frac{ \|\mathcal C[G;\mathbb R]^{-1}\|^2}{1 - \epsilon \|\mathcal C[G;\mathbb R]^{-1}\|}.
\end{align*}
Furthermore, if $\mathcal C[G;\mathbb R] u = G-I+F$ and $\mathcal C[\hat G; \mathbb R]\hat u = \hat G -I + \hat F$ then
\begin{align*}
\|u - \hat u\|_{L^2(\mathbb R)} &\leq \epsilon \|\hat F\|_{L^2(\mathbb R)} \frac{ \|\mathcal C[G;\mathbb R]^{-1}\|^2}{1 - \epsilon \|\mathcal C[G;\mathbb R]^{-1}\|} + \epsilon \|\mathcal C[G;\mathbb R]\| \triangleq B(\epsilon).
\end{align*}
If, in addition,
\begin{align*}
\left|\dashint_{\mathbb R} (\hat G(x) -G(x))dx \right| < \epsilon \quad \text{and} \quad  \left|\dashint_{\mathbb R} (\hat F(x) -F(x))dx \right| < \epsilon,
\end{align*}
then
\begin{align*}
\left|\dashint_{\mathbb R} (\hat u(x) -u(x)) dx \right| &\leq 2 \epsilon + \epsilon \|u\|_{L^2(\mathbb R)} + B(\epsilon) \|\hat G -I\|_{L^2(\mathbb R)}.
\end{align*}
\begin{proof}
All statements aside from the last follow from the application of standard results in operator theory \cite{atkinson} (see also \cite{TrogdonThesis} for the case of singular integral equations on general contours).  The final statement follows from an application of the Cauchy--Schwarz inequality when considering the difference $u-\hat u$ using the equations that $u$ and $\hat u$ satisfy.
\end{proof}
\end{lemma}

If $\hat G$ is a finite sum of the basis $\{R_j\}$ the methods described above allow for the application of $\mathcal C[\hat G;\mathbb R]$ exactly to a function that is a finite sum of the basis $\{R_j\}$ and this process returns a function that is again a finite sum of the basis $\{R_j\}$.  Additionally, inner products of such functions are computed exactly and hence the infinite-dimensional GMRES algorithm may be applied \cite{GMRES-original}.

\section{Infinite-dimensional GMRES}\label{Section:GMRES}
\def\ipH<#1><#2>{\langle #1, #2 \rangle_*}

Let $\mathbb H$ be a Hilbert space with inner product $\ipH<\cdot><\cdot>$.  The fundamental idea of GMRES for the solution of an operator equation $\mathcal A x = f$ with $x,f \in \mathbb H$, is the solution of the minimization problem
\begin{align*}
\inf_{x \in \mathbb K_n} \|A x -f\|_*, ~~ \mathbb K_n = \text{span} \{f,\mathcal Af,\mathcal A^2f,\ldots,\mathcal A^{n-1}f\}.
\end{align*}
The solution of this problem gives $x_n$ and ideally $\|x_n - x^*\|_* \goto 0$ where $x^*$ is the true solution.

\subsection{Arnoldi Iteration}

For the stable solution of the above minimization problem we use the Arnoldi algorithm.  The algorithm expresses the action of the operator $\mathcal A$ on $\mathbb K_n$ in the form
\begin{align*}
\mathcal A Q_n = Q_{n+1} \tilde{H_n}.
\end{align*}
The columns (elements of $L^2(\mathbb R)$) of $Q_n$ form an orthonormal basis for $\mathbb K_n$ and
\begin{align*}
\tilde{ H_n} = \begin{mat} h_{11} & h_{12} & h_{13}&\cdots  & h_{1,n} \\
h_{21} & h_{22}  &&&\vdots\\
&\ddots&\ddots &&\vdots\\
&  & \ddots &\ddots&\vdots\\
&&& h_{n,n-1} & h_{n,n} \end{mat}
\end{align*}
is upper Hessenberg.

\subsection{Mechanics of GMRES}

Again, consider the minimization of
\begin{align*}
\|\mathcal A x - f\|^2_*, ~~ x \in \mathbb K_n.
\end{align*}
By expressing $x = Q_n y$ for $y \in \mathbb C^n$ and $f = \|f\|_* Q_{n+1} e_1$ where $e_1 = (1,0,\ldots,0)^T$ one is led to the minimization of
\begin{align*}
\|\mathcal A Q_n y - \|f\|_* Q_{n+1} e_1\|^2_*, ~~ y \in \mathbb C^n.
\end{align*}
This is further reduced to
\begin{align*}
\|Q_{n+1} (\tilde H_n y - \|f\|_* e_1)\|^2_*, ~~ y \in \mathbb C^n.
\end{align*}
Since $Q_{n+1}$ is an isometry from $\mathbb C^{n+1}$ (with the usual $l^2$ inner product) to $\mathbb K_{n+1}$, this is equivalent to the minimization of
\begin{align*}
\|\tilde H_n y-  \|f\|_* e_1\|_2^2, ~~~ y \in \mathbb C^n.
\end{align*}
Once we have $y$, the approximate solution $x_n$ of $\mathcal Ax=f$ is $x_n = Q_n y$.

We need the $QR$ factorization of $\tilde H_n$ for each $n$ to solve the minimization problem.  There exists multiple ways of computing this factorization.  Since $\tilde H_n$ is being built iteratively,  Givens rotations and $\bigo(n)$ work build successive $QR$ factorizations. We see that few iterations of GMRES are needed so that the method used to solve the least squares problem is of lower importance.

\subsection{Convergence}

We give a brief discussion of what is known about the application of GMRES in infinite dimensions.  We say an operator $\mathcal A$ is an algebraic operator if there exists a polynomial $p(x) \in \mathbb C[x]$ such that $p(\mathcal A) = 0$.  A sufficient condition for this is if $\mathcal A = \lambda \id + \mathcal K$ where $\lambda \in \mathbb C$ and $\mathcal K$ is of finite rank \cite{Brevsar}.  GMRES applied (exactly) to $\mathcal Ax=f$ in an infinite-dimensional separable Hilbert space is known to converge if $\mathcal A$ is an algebraic operator \cite{Gasparo}.  In the examples below we apply GMRES in cases when $A$ is clearly algebraic and when it is not known if $A$ is algebraic.  Convergence is demonstrated on a case-by-case basis.


\section{Non-oscillatory examples and results}\label{Section:NonOscExamples}

We pause briefly to discuss some implementational details of the method.  With each application of $\mathcal C[\hat G;\mathbb R]$ ($\hat G$ is expressed in terms of the basis $\{R_j\}_{j=-N}^N$) we obtain a function that is expressed in terms of $2N$ more basis functions due to function multiplication.  In principle, this fact can be ignored.  In practice, truncation should be employed to increase efficiency.  With each application of the operator, all coefficients are dropped that fall below a given tolerance $\epsilon_{\text{trunc}}>0$.  This process keeps the computational cost under control.  We fix a tolerance $\epsilon> 0$ and halt GMRES when the residual falls below $\epsilon$.  We always take $\epsilon > \epsilon_{\text{trunc}}$.  If this is not the case, truncation may cause GMRES to fail to converge.

\subsection{A scalar problem}\label{Section:scalar}

\begin{figure}[tb]
\centering
\subfigure[]{\includegraphics[width=.45\linewidth]{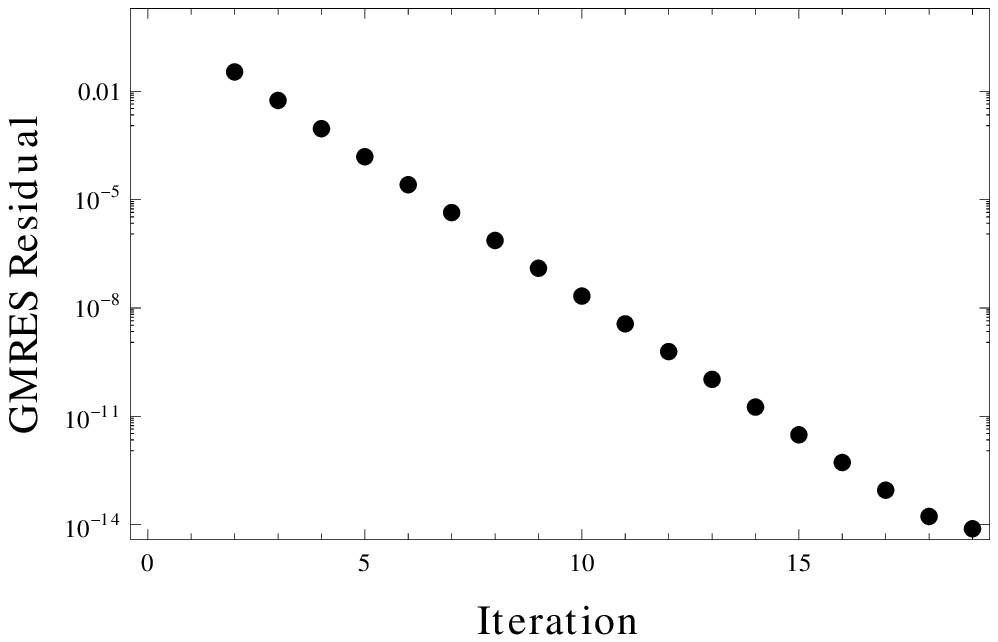}\label{ScalarGMRES-naive}}
\subfigure[]{\includegraphics[width=.45\linewidth]{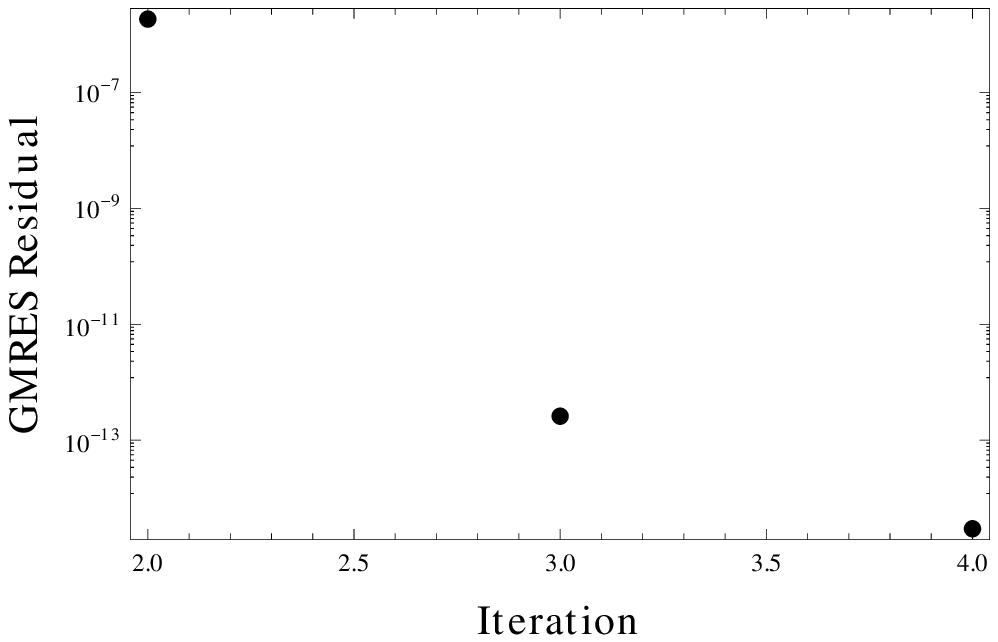}\label{ScalarGMRES-reg}}
\caption{(a) Convergence of the GMRES residual for $\mathcal C[1 + \sech(\cdot);\mathbb R]u(x) = \sech(x)$.  (b) Convergence of the GMRES residual for \eqref{scalar-precond}.  Only four iterations are needed for GMRES to converge to machine precision when the preconditioner is used.}
\end{figure}

Consider the scalar RHP
\begin{align}\label{scalar-rhp}
\Phi^+(x) = \Phi^-(x) (1+\sech(x)), ~~ \Phi(\infty) = 1.
\end{align}
From \eqref{sie} we consider the singular integral equation
\begin{align*}
\mathcal C[1+\sech(\cdot);\mathbb R]u(x) = u(x)- \mathcal C_{\mathbb R}^-u(x) \cdot\sech(x) = \sech(x).
\end{align*}
We approximate $\sech(x)$ in the basis $\{R_j\}_{j=-N}^N$ with $N = 250$.  See Figure~\ref{ScalarGMRES-naive} for the convergence of the GMRES residual.

The theory of singular integral equations suggests a preconditioner.  The operator
\begin{align*}
\mathcal C[1/(\sech(\cdot ) + 1);\mathbb R]u(x)= u(x) = \mathcal C_{\mathbb R}^- u(x) (1-1/(\sech(x)+1))
\end{align*}
is a Fredholm regulator for $\mathcal C[1 + \sech(\cdot);\mathbb R]$ in the sense that 
\begin{align*}
\mathcal C[1/(\sech(\cdot) + 1);\mathbb R]\mathcal C[1 + \sech(\cdot);\mathbb R] = \id + \mathcal K
\end{align*}
where $\mathcal K$ is a compact operator on $L^2(\mathbb R)$ \cite{TrogdonThesis,zhou-RHP}.  Furthermore,  if we replace $\sech(x)$ and $1/(1+\sech(x))-1$ with their rational approximations in the basis $\{R_j\}$, $\mathcal K$ is a finite-rank operator \cite{TrogdonThesis}.  See Figure~\ref{ScalarGMRES-reg} for the convergence of the GMRES residual associated with the equation
\begin{align}\label{scalar-precond}
\mathcal C[1/(\sech(\cdot) + 1);\mathbb R]\mathcal C[1 + \sech(\cdot);\mathbb R]u(x) = \mathcal C[1/(\sech(\cdot) + 1);\mathbb R]\sech(x).
\end{align}
 We see that only four iterations of GMRES are needed in this case.  This presents an important link between the analysis of singular integral operators and numerical analysis.  This link is further emphasized in future examples.

\begin{figure}[tb]
\centering
\includegraphics[width=.5\linewidth]{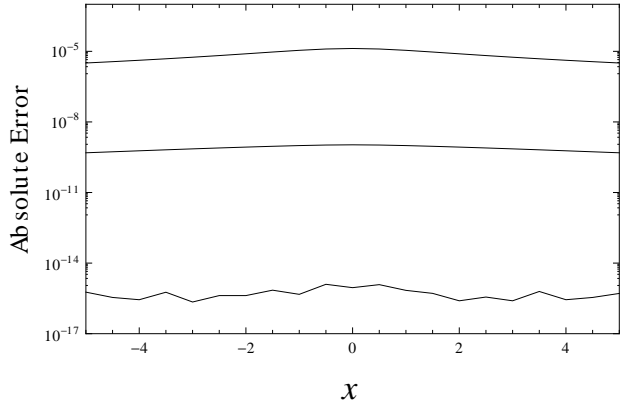}
\caption{\label{ScalarCompare} Convergence of $\mathcal C_{\mathbb R}^+u_n$ to the true solution of \eqref{scalar-rhp} where $n$ represents the number of iterations of GMRES that have been performed.  These results are in the absence of the preconditioner.  Absolute error is plotted versus $x$ for $n=5,10,19$.  }
\end{figure}

This RHP can be solved explicitly
\begin{align*}
\Phi(z) = \exp\left( \frac{1}{2\pi i} \int_{\mathbb R} \frac{\log(1 + \sech(s))}{s-z} ds   \right),
\end{align*}
and $1+\mathcal C^\pm_{\mathbb R}u = \Phi^\pm$.  This expression can be evaluated accurately using the method described above for the Cauchy integral to provide a comparison for the method here.  See Figure~\ref{ScalarCompare} for a demonstration of the convergence of GMRES to this solution.

\subsection{A matrix problem}

We move to consider matrix singular integral equations.  Again, we concentrate on equations that arise in the solution of RHPs.  Consider the $2\times 2$ matrix RHP
\begin{align}\label{IST-RHP}
\begin{split}
\Phi^+(z) = \Phi^-(z) G(z;x,t), ~~ z \in \mathbb R, ~~\Phi(\infty) = I,\\
G(z;x,t) = \begin{mat} 1 - |\rho(z)|^2 & - \bar\rho(z) e^{-2ixz-4iz^2t},\\
\rho(z) e^{2ixz+4iz^2t} & 1 \end{mat}.
\end{split}
\end{align}
No explicit solution of this problem is known.  This RHP arises in the solution of the defocusing nonlinear Schr\"odinger (NLS) equation\begin{align}\label{NLS}
\begin{split}
-iq_t+q_{xx} - 2q|q|^2 &= 0,\\
q(x,0) &= q_0(x),
\end{split}
\end{align} 
with the inverse scattering transform \cite{AblowitzClarksonSolitons}.  Here $\rho$ is the reflection coefficient associated with a decaying initial condition $q_0$.  See \cite{TrogdonSONLS} for a discussion of the computation of $\rho$ given an initial condition $q_0$.  Once this RHP is solved the solution of the NLS equation is found through the formula
\begin{align*}
q(x,t) = -2i\lim_{z\goto \infty} z \Phi_{21}(z),
\end{align*}
where the subscripts denote the $(2,1)$ entry of $\Phi$.

\begin{figure}[tb]
\centering
\subfigure[]{\includegraphics[width=.45\linewidth]{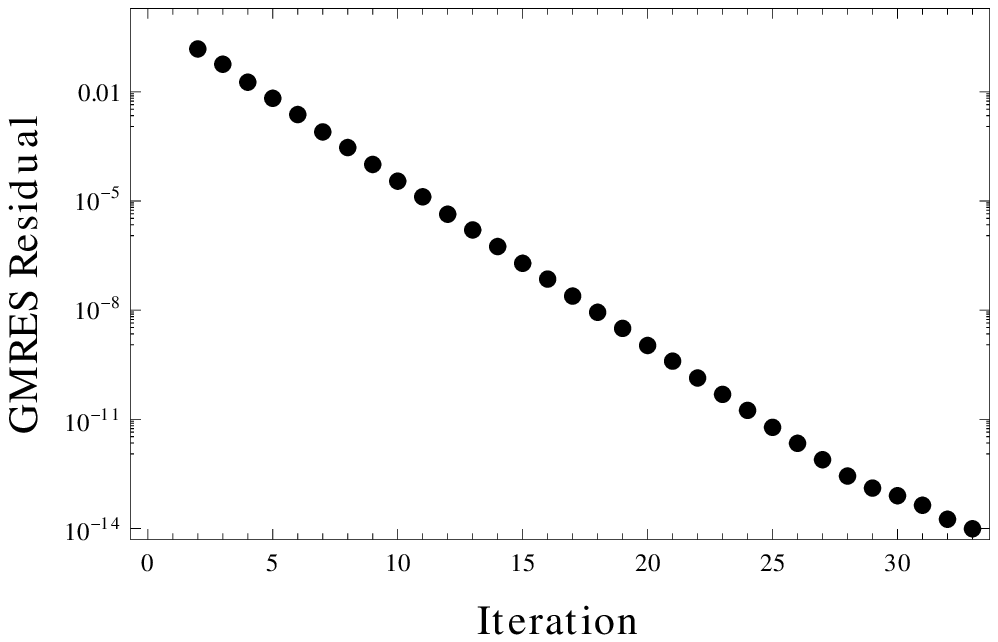}\label{MatrixGMRES-naive}}
\subfigure[]{\includegraphics[width=.45\linewidth]{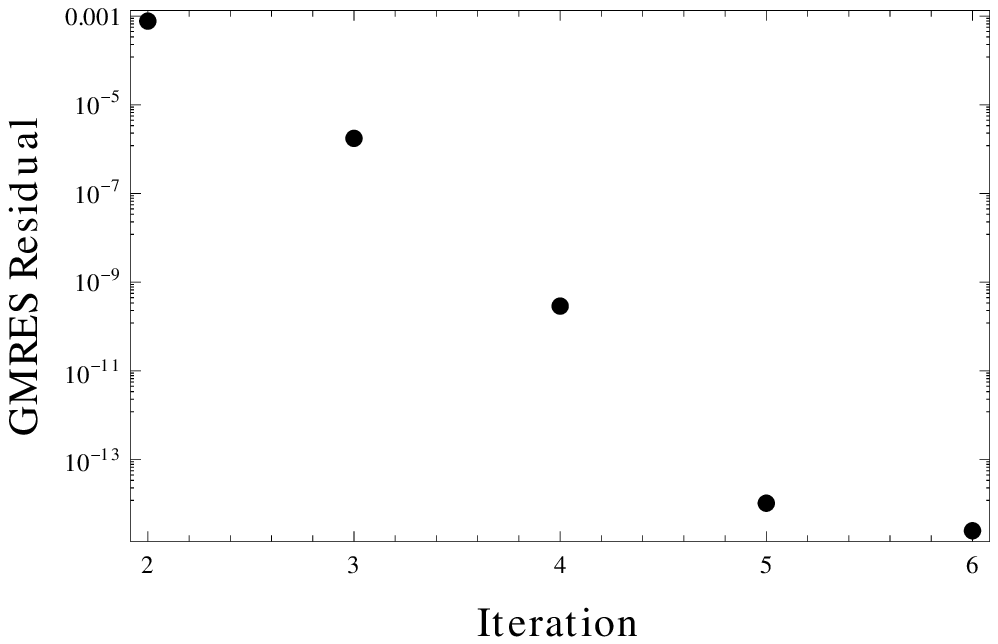}\label{MatrixGMRES-reg}}
\caption{(a) Convergence of the GMRES residual for $\mathcal C[G;\mathbb R]u = G-I$. (b) Convergence of the GMRES residual for $\mathcal C[G^{-1};\mathbb R]\mathcal C[G;\mathbb R]u = \mathcal C[G^{-1};\mathbb R](G-I)$.  Only six iterations are need to achieve machine accuracy for the preconditioned operator.}
\end{figure}

It is known that if $q_0$ is smooth and rapidly decaying (faster than any polynomial) then so is $\rho$ \cite{deift-zhou:nls}.  As a prototypical example we use $q_0(x) = e^{-x^2}$.  For small values of $|x|$ and $|t|$, $\rho(z)e^{2ixz+4iz^2t}$ may be accurately expressed in terms of the basis $\{R_j\}$ and we may apply GMRES to the operator equation
\begin{align*}
\mathcal C[G;\mathbb R]u = G-I,
\end{align*}
to solve the RHP.  See Figure~\ref{MatrixGMRES-naive} for the convergence of the GMRES residual.  We may also apply GMRES to the equation
\begin{align*}
\mathcal C[G^{-1};\mathbb R] \mathcal C[G;\mathbb R]u = \mathcal C[G^{-1};\mathbb R](G-I),
\end{align*}
and faster convergence is realized.  In what follows we refer to this as the preconditioned equation. Indeed, like the scalar case  $\mathcal C[G^{-1};\mathbb R] \mathcal C[G;\mathbb R]=\id + \mathcal K$ where $\mathcal K$ is compact \cite{TrogdonThesis}. See Figure~\ref{MatrixGMRES-reg} for the convergence of the GMRES residual for the preconditioned equation.

\begin{figure}[tb]
\centering
\includegraphics[width=.45\linewidth]{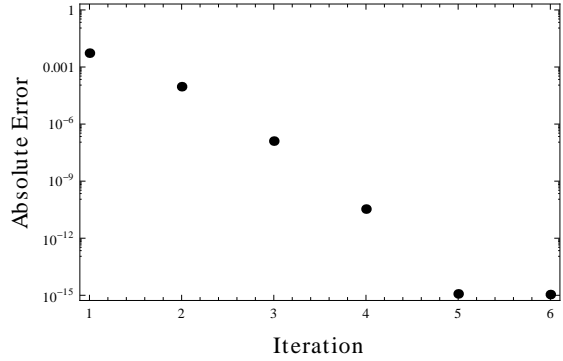}\label{zero-zero-error}
\caption{A demonstration of the error in the approximation of $q(0,0)$ at each iteration of GMRES.}
\end{figure}

It follows that
\begin{align}\label{integral-reconstruct}
q(x,t) = -2i\lim_{z \goto \infty} z \Phi_{21} = \frac{1}{\pi} \int_{\mathbb R} u_{21}(z) dz,
\end{align}
and the methods presented above allow us to approximate such an integral when we have an approximation of $u$ in terms of the basis $\{R_j\}$.  To demonstrate the convergence of the method we set $x=t=0$ and compute an approximation of $q(0,0)$ using \eqref{integral-reconstruct} for each iteration of GMRES.  Numerical results are shown in Figure~\ref{zero-zero-error} in the case of the preconditioned compact operator.  Accuracy on the order of machine precision is easily obtained.

\begin{remark}\label{non-osc-issues}
For $ t >0$ small and $|x|$ small, we might approximate $e^{2izx+4iz^2t}\rho(z)$ with $\{R_j\}$ and solve the singular integral equation \eqref{sie}. Two complications are present. First, as $|x|$ and $t$ increase more basis functions are needed to resolve the jump matrix and hence, more are needed to resolve the solution. This slows the computation significantly.  Secondly, the operator tends to be ill-conditioned in the sense that it takes GMRES many iterations to converge. In practice, we are limited to $|x|\leq 2$ and $t<0.01$ due to time constraints.  This restriction on both $|x|$ and $t$ being small is important. In the following sections we remove the restriction on $|x|$ small by introducing oscillatory basis functions.  It is still unknown how to deal with time dependence in a similar way due to the quadratic nature of the oscillations..
\end{remark}

\section{Oscillatory basis functions}\label{Section:OscBasis}

As the examples demonstrate, when the non-oscillatory method above is applied to the singular integral equations associated with inverse scattering the acceptable range of parameter values ($x$ and $t$) is limited.  By introducing a two-parameter family of oscillatory basis functions, we make progress on removing this restriction.  Define (for fixed $\beta > 0$)
\begin{align*}
R_{j,\alpha}(z) = e^{i \alpha z}\left( \left(\frac{z-i\beta}{z+i\beta}\right)^j - 1\right),  ~~ \alpha \in \mathbb R.
\end{align*}
It is clear that the set
\begin{align}\label{osc-basis}
\{R_{j,\alpha}(z) \}_{j \in \mathbb Z, ~\alpha \in \mathbb R},
\end{align}
is a generalization of the non-oscillatory basis \eqref{non-osc-basis}.  Note that this basis is convenient for representing $\rho(z)e^{2izx + 4iz^2t}$ for $t = 0$, $\alpha = 2x$.  It is important that we set $t = 0$ so that there is no quadratic term in the phase.   We refer to \eqref{osc-basis} as the oscillatory basis.  As before, the action of the Cauchy operators on this basis, how multiplication transforms the basis and the computation inner products must be understood. We mirror Section~\ref{Section:non-osc}.

\subsection{Properties of the oscillatory basis}

First, it is clear that the oscillatory basis is dense in $L^2(\mathbb R)$ since it contains \eqref{non-osc-basis}.  Multiplication is also straightforward using \eqref{non-osc-mult}:
\begin{align*}
R_{j,\alpha_1}(z)R_{k,\alpha_2}(z) = R_{k+j,\alpha_1+\alpha_2}(z) -R_{j,\alpha_1+\alpha_2}(z)-R_{k,\alpha_1 + \alpha_2}(z).
\end{align*}
We concentrate on the other properties.

\subsubsection{Action of the Cauchy operators}

We compute the action of the Cauchy operators through residue calculations.  The following lemma assists in these calculations.

\begin{lemma}\label{Lemma:Residue}
For $z$ sufficiently close to the real axis
\begin{align*}
\res\left\{ R_{j,\alpha}(s)\frac{1}{s-z}; s = z\right\} = R_{j,\alpha}(z),
\end{align*}
and for $\sigma = \sign(j)$
\begin{align*}
\res&\left\{ R_{j,\alpha}(s)\frac{1}{s-z}; s = -\sigma i\beta\right\} = \sum_{n=0}^{|j|} \gamma_n\frac{(-2 i \sigma \beta)^n}{(z+\sigma i \beta)^{n+1}},\\
\gamma^j_n& = - e^{\sigma \alpha \beta} \sum_{k=0}^{|j|-n}\left(\begin{array}{c} |j|-n \\ k \end{array} \right) \left(\begin{array}{c} |j| \\ |j|-n \end{array} \right) \frac{n!}{(n+k)!} (2\sigma \beta \alpha)^k,\\
\res&\left\{ R_{j,\alpha}(s)\frac{1}{s-z}; s = \sigma i\beta\right\} = 0.
\end{align*}
\begin{proof}
We prove the result for $j >0$ as $j < 0$ requires only the addition of absolute value signs.  The only non-trivial statement is the calculation when $s=-\sigma i\beta$.  We must find the $(j-1)$th term in the Taylor expansion of $e^{i\alpha s}(s-\sigma i\beta)^j(s-z)^{-1}$ about $s = -\sigma i\beta$. The $\ell$th derivative of the first two factors of this function is
\begin{align}\label{taylor}
\frac{d^\ell}{ds^\ell} \left. \left( e^{i\alpha s}(s-\sigma i\beta)^j \right)\right|_{s = -\sigma i\beta} = \sum_{k=0}^\ell \left(\begin{array}{c} \ell \\ k \end{array} \right) \frac{j!}{(j-(\ell-k))!} (i \alpha)^k (-2i\sigma\beta)^{j-(\ell-k)} e^{\sigma \alpha \beta}.
\end{align} 
We obtain a double sum for the $(j-1)$th derivative
\begin{align*}
d_{j-1}&\triangleq \frac{d^{j-1}}{ds^{j-1}} \left. \left( e^{i\alpha s}(s-\sigma i\beta)^j \frac{1}{s-z} \right)\right|_{s = -\sigma i\beta} \\
&= -e^{\sigma \alpha \beta}\sum_{\ell = 0}^{j-1} \sum_{k=0}^\ell \left(\begin{array}{c} \ell \\ k \end{array} \right) \left(\begin{array}{c} j-1 \\ \ell \end{array} \right) \frac{j!(j-\ell-1)!}{(j-\ell+k)!} \frac{(i\alpha)^k (-2 \sigma i \beta)^{j-\ell+k}}{(z+\sigma i\beta)^{j-\ell}}.
\end{align*}
We perform the change of variables $n = j-\ell$ and obtain
\begin{align*}
\frac{d_{j-1}}{(j-1)!} &=   \sum_{n=1}^j \gamma^j_n\frac{(-2 i \sigma \beta)^n}{(z+\sigma i \beta)^{n}},\\
\gamma^j_n& = - j e^{\sigma \alpha \beta} \sum_{k=0}^{j-n}\left(\begin{array}{c} j-n \\ k \end{array} \right) \left(\begin{array}{c} j-1 \\ j-n \end{array} \right) \frac{(n-1)!}{(n+k)!} (2\sigma \beta \alpha)^k.
\end{align*}
\end{proof}
\end{lemma}

We make use of another lemma for the relationship between bases.  
\begin{lemma}\label{Lemma:Basis}
For $\sigma = \pm 1$ and $j > 0$ we have
\begin{align}\label{to-simple}
R_{\sigma j,\alpha}(z) = \sum_{n=1}^{j} \left(\begin{array}{c} j \\ n \end{array} \right) \frac{(-2i\sigma \beta)^{n}}{(z+\sigma i\beta)^{n}} e^{i \alpha z},
\end{align}
and
\begin{align*}
\frac{(-2i\sigma \beta)^{j}}{(z+\sigma i\beta)^{j}} e^{i \alpha z} = \sum_{n=1}^{j} (-1)^{j+n}\left(\begin{array}{c} j \\ n \end{array} \right) R_{\sigma n,\alpha}(z).
\end{align*}
\begin{proof}
We write
\begin{align*}
R_{j,\alpha}(z) = \frac{(z-i\beta)^j - (z+i\beta)^j}{(z+i\beta)^j}.
\end{align*}
Using
\begin{align*}
(z+i\beta)^j = \sum_{\ell = 0}^j \left(\begin{array}{c} j \\ \ell \end{array} \right) (z+i\beta)^\ell (-2i\beta)^{j-\ell},
\end{align*}
we find \eqref{to-simple}.  Next, the inverse of the matrix
\begin{align*}
A_{ij} = \begin{choices} \left(\begin{array}{c} j \\ n \end{array} \right), \when i \leq j,\\
0, \otherwise, \end{choices}
\end{align*}
is given by
\begin{align*}
A^{-1}_{ij} = \begin{choices} (-1)^{j+n}\left(\begin{array}{c} j \\ n \end{array} \right), \when i \leq j,\\
0, \otherwise, \end{choices}
\end{align*}
\end{proof}
\end{lemma}

For efficient computation, we simplify the expression for $\gamma^j_n$ using Krummer's confluent hypergeometric function \cite{DLMF}
\begin{align}\label{1F1}
_1F_1(a,b,z) = \sum_{k=0}^\infty \frac{(a)_k}{(b)_k} \frac{z^k}{k!},
\end{align}
where $(a)_k$ is the Pochhammer symbol \cite{DLMF},
\begin{align*}
(a)_k = \prod_{i=0}^{k-1} (a+i) = \frac{\Gamma(a+k)}{\Gamma(a)},
\end{align*}
and $\Gamma(z)$ is the Gamma function.  Note that if $a < 0$ and $k \geq a$ then $(a)_k = 0$ and this truncates \eqref{1F1} to a finite sum.  Properties of the Gamma function can be used to show that
\begin{align*}
\frac{\Gamma(n-j+k)\Gamma(j-n-k+1)}{\Gamma(n-j)\Gamma(j-n+1)} = (-1)^k.
\end{align*}
From this it follows that
\begin{align*}
\gamma^j_n = - \frac{|j|}{n} e^{\sigma \alpha \beta} \left(\begin{array}{c} |j|-1 \\ n \end{array} \right) \phantom{.}_1F_1(n-|j|,1+n,-2\sigma \alpha \beta).
\end{align*}
We arrive at the following lemma that shows how taking a residue maps the basis to itself.
\begin{lemma}\label{Lemma:eta}
For $\sigma = \sign(j)$
\begin{align}\label{res-eta}
\res&\left\{ R_{j,\alpha}(s)\frac{1}{s-z}; s = -\sigma i\beta\right\} = \sum_{n=1}^{|j|} \eta^j_n R_{\sigma n,0}(z),\notag\\
\eta^j_n &= \sum_{k=n}^{|j|}(-1)^{n+k} \left(\begin{array}{c} k \\ n \end{array} \right)  \gamma_k^j.
\end{align}
\begin{proof}
This follows directly from Lemmas \ref{Lemma:Residue} and \ref{Lemma:Basis}.
\end{proof}

\end{lemma}

\begin{remark}
It is clear that $\eta^j_n$ depends on $\alpha$, $\beta$ and $\sigma = \sign(j)$ but we suppress these parameters for ease of notation.
\end{remark}

\begin{remark}
The series expression for $\eta_n^j$ is an alternating series.  Indeed, $\gamma_n^j$ is also an alternating series.  Therefore it is difficult to re-order the sum in such a way to explicitly sum the alternating terms first.  Stable computation of $\eta_n^j$ is difficult.  In practice, we use higher-precision arithmetic to compute $_1F_1$ accurately and accurately compute $\eta_n^j$.  We treat the computation of $\eta_n^j$ as black-box special function. 
\end{remark}

With this lemma in hand we are able to describe the action of the Cauchy operators on the oscillatory basis.

\begin{theorem}
If $\alpha \cdot j \geq 0$ then
\begin{align*}
\mathcal C^+_{\mathbb R} R_{j,\alpha}(z) &= \begin{choices} R_{j,\alpha}(z), \when j \geq 0,\\ 0, \when j < 0, \end{choices}\\
\mathcal C^-_{\mathbb R} R_{j,\alpha}(z) &= \begin{choices} 0, \when j \geq 0, \\ -R_{j,\alpha}(z), \when j < 0.\end{choices}
\end{align*}
If $\alpha \cdot j < 0$ then
\begin{align*}
\mathcal C^+_{\mathbb R} R_{j,\alpha}(z) & = \begin{choices} - \displaystyle \sum_{n=1}^{j} \eta^j_n R_{n,0}(z), \when j \geq 0,\\ 
R_{j,\alpha}(z) +  \displaystyle\sum_{n=1}^{-j} \eta^j_n R_{-n,0}(z), \when j < 0, \end{choices}\\
\mathcal C^-_{\mathbb R} R_{j,\alpha}(z) & = \begin{choices} -R_{j,\alpha}(z)-  \displaystyle\sum_{n=1}^{j} \eta^j_n R_{n,0}(z), \when j \geq 0,\\ 
  \displaystyle\sum_{n=1}^{-j} \eta^j_n R_{-n,0}(z), \when j < 0. \end{choices}
\end{align*}
\begin{proof}
This follows from straightforward residue calculations using Lemmas~\ref{Lemma:Residue} and \ref{Lemma:eta} and Cauchy's Theorem.
\end{proof}
\end{theorem}

\begin{remark}  Note that the coefficients $\eta_n^j$ appear only when $\alpha \sign(j)  < 0$.  This is sufficient to ensure that the exponential in \eqref{res-eta} always induces decay.
\end{remark}

\subsubsection{The inner product}

Again, looking toward the application of GMRES we compute the inner products
\begin{align*}
L_{j,k,\alpha_1,\alpha_2} = \int_{\mathbb R} R_{j,\alpha_1}(z) \overline{R_{k,\alpha_2}(z)} dz.
\end{align*}
Since $\overline{R_{k,\alpha_2}(z)} = R_{-k,-\alpha_2}(z)$
\begin{align*}
L_{j,k,\alpha_1,\alpha_2} = \dashint_{\mathbb R} (R_{j-k,\alpha_1-\alpha_2}(z) - R_{j,\alpha_1-\alpha_2}(z) - R_{-k,\alpha_1-\alpha_2}(z))dz.
\end{align*}
Therefore, the problem reduces to computing
\begin{align}\label{osc-pv}
\dashint_{\mathbb R} R_{j,\alpha}(z) dz. 
\end{align}
This is nothing more than the Fourier transform of the non-oscillatory basis. Jordan's Lemma \cite[p.~222]{FokasComplexVariables} shows us that
\begin{align*}
\int_{C_r^+} R_{j,\alpha}(z) dz= 0, ~~ \alpha > 0, ~~ \int_{C_r^-} R_{j,\alpha}(z) dz=0, ~~ \alpha < 0,
\end{align*}
where $C_r^\pm$ is defined below \eqref{c-r}.  Computing \eqref{osc-pv} reduces to a pure residue calculation.  Using \eqref{taylor} for $\sigma = \pm 1$, $ j >0$ and $\alpha \neq 0$
\begin{align*}
\dashint_{\mathbb R} R_{\sigma j,\alpha}(z) dz &= 2 \pi i \sign(\alpha) \res\{ R_{\sigma j,\alpha}(z), z = \sign(\alpha) i \beta \}\\
&= \begin{choices}
0, \when  \sign(\alpha) = \sigma,\\
\displaystyle \sign(\alpha) \sum_{k=0}^{j-1} \left(\begin{array}{c} j-1 \\ k \end{array} \right) \frac{j!}{(k+1)!} (i \alpha)^k (2i\sign(\alpha)\beta)^{k+1} e^{-|\alpha| \beta}, \otherwise.\end{choices}
\end{align*}
This sum can be simplified using the $_1F_1$ function.  We find
\begin{align}\label{osc-formula}
I_{j,\alpha} \triangleq \dashint_{\mathbb R} R_{j,\alpha}(z) dz = \begin{choices}
0, \when \sign(j) = \sign(\alpha),\\
- 2 \pi |j| \beta, \when \alpha = 0,\\
-4 \pi e^{- |\alpha| \beta} |j| \beta \phantom{.}_1F_1(1-|j|,2,2 |\alpha| \beta ), \otherwise.\end{choices}
\end{align}
Therefore
\begin{align}\label{osc-ip}
L_{j,k,\alpha_1,\alpha_2} = I_{j-k,\alpha_1-\alpha_2} - I_{j,\alpha_1-\alpha_2} - I_{-k,\alpha_1-\alpha_2}.
\end{align}

\begin{remark}
Since $\{R_j(z)\}_{j=-\infty}^\infty$ is a basis of $L^2(\mathbb R,dz)$ and the Fourier transform is unitary on $L^2(\mathbb R,dz)$ we have that $\{\chi_{\{\alpha j <0\}}(\alpha) e^{- |\alpha| \beta} \phantom{.}_1F_1(1-j,2,2 |\alpha| \beta )\}_{j=-\infty}^{\infty}$ is a basis of $L^2(\mathbb R, d\alpha)$ where $\chi_A$ is the characteristic function of the set $A$.  In Figure~\ref{new-basis} we plot basis functions for a couple values of $j$.
\begin{figure}[tb]
\centering
\subfigure[]{\includegraphics[width=.45\linewidth]{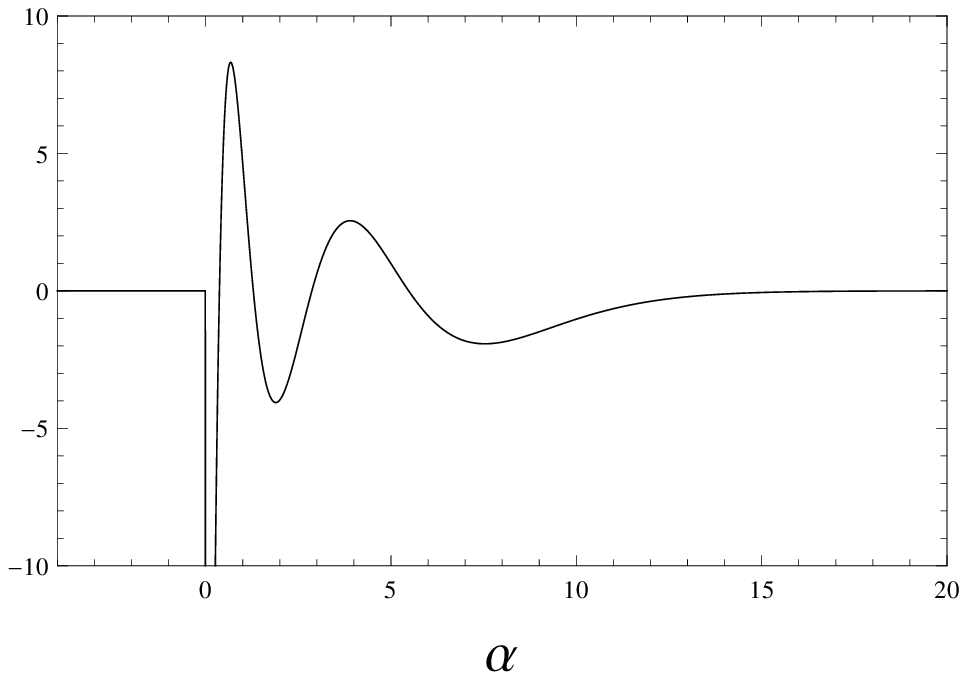}}
\subfigure[]{\includegraphics[width=.45\linewidth]{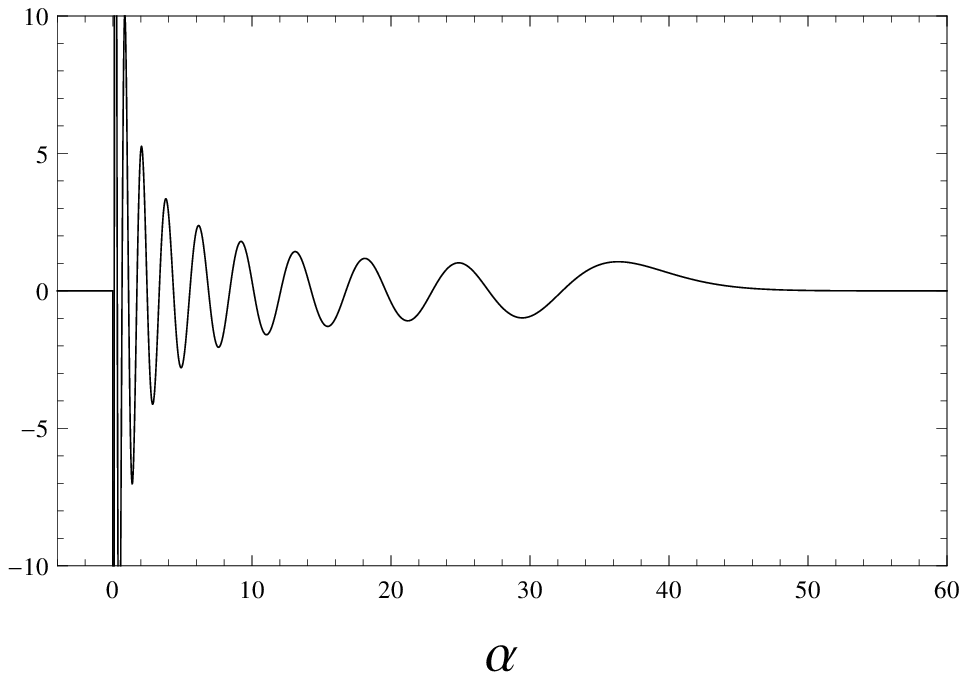}}
\caption{\label{new-basis} (a) A plot of the Fourier transform of $R_{-5}(z)$.  (b) A plot of the Fourier transform of $R_{-20}(z)$.}
\end{figure}
\end{remark}

\begin{remark}
An alternate representation of the basis can be derived because \cite{Abramowitz}
\begin{align*}
_1F_1(-n,a+1,z) = \frac{n!}{(a+1)_n} L^{(a)}_n(x), 
\end{align*}
where $L^{(a)}_n(x)$ is the generalized Laguerre polynomial of order $n$.  In practice, we see that the stock methods in {\tt Mathematica} are accurate for $ L^{(a)}_n(x)$ for large $n$ but stock methods show inaccuracies for $_1F_1$ for large negative values of the first parameter.
\end{remark}

\subsubsection{Application to oscillatory quadrature}

It is clear from the previous section that \eqref{osc-formula} has application to numerical Fourier analysis.  Here we present two examples to demonstrate these applications.  An analysis of this method can be found in \cite{Trogdon2014}.

\begin{example}[Computing Fourier transforms]
We use \eqref{osc-formula} to compute the Fourier transform of a Gaussian $f(x) = e^{-x^2}$ with $\beta = 1$.  We follow Example~\ref{Example:RationalApprox} to approximate $f$ with the basis $\{R_j\}$. This produces the coefficients of the Fourier transform in the basis $\{\chi_{\{\alpha j <0\}}(\alpha) e^{- |\alpha| \beta} \phantom{.}_1F_1(1-j,2,2 |\alpha| \beta )\}_{j=-\infty}^{\infty}$.  See Figure~\ref{FT-error} for a demonstration of the accuracy of the method.

\begin{figure}[btp]
\centering
\includegraphics[width=.5\linewidth]{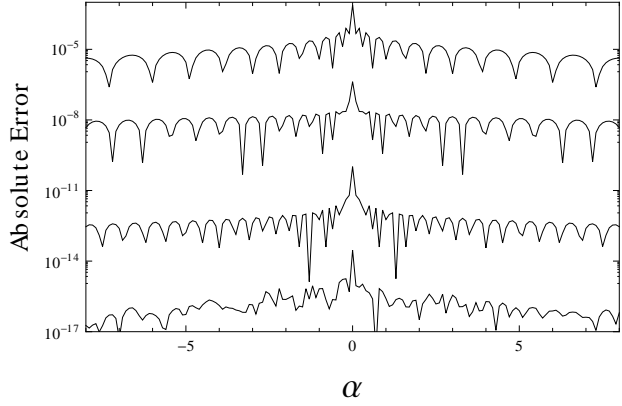}
\caption{\label{FT-error} Error in the computation of the Fourier transform of $f(x)=e^{-x^2}$ when $f$ is expanded in the basis $\{R_j\}_{j=-N}^N$.  Absolute error is plotted versus $\alpha$ for $N=20,40,80$ and $160$.  Compare this with Figure~\ref{Gaussian-approximation} to see that there is no loss of accuracy for $|\alpha| > 1$. }
\end{figure}
\end{example}

\begin{example}
We also use these techniques to solve linear partial differential equations for small time.  Consider the PDE
\begin{align}\label{LS}
-iq_t + q_{xx} = 0,\\
q_0(x) = e^{-x^2}\notag,
\end{align}
posed on $\mathbb R$.  We solve this PDE with the Fourier transform and explicitly compute the transform of the initial condition:
\begin{align*}
\hat q_0(z) = \int_{\mathbb R} e^{-izx}q_0(x) dx = \sqrt{\pi} e^{-z^2/4}.
\end{align*}
We approximate $\hat q_0(z) e^{iz^2t}$ with the basis $\{R_j\}_{j=-\infty}^\infty$ with $\beta = 1$:
\begin{align*}
\hat q_0(z) e^{iz^2t} \approx \sum_{j=-N}^N \gamma_j R_j(z).
\end{align*}
This approximation is only viable for small $t$.  The techniques described above allow us to compute the approximation
\begin{align*}
q(x,t) \approx \frac{1}{2 \pi} \sum_{j=-N}^N \gamma_j \dashint_{\mathbb R} R_{j,x}(z) dz
\end{align*}
which is uniformly valid in $x$.  See Figure~\ref{FT-soln} for a plot of the solution of \eqref{LS} computed with this method.

\begin{figure}[tbp]
\centering
\subfigure[]{\includegraphics[width=.45\linewidth]{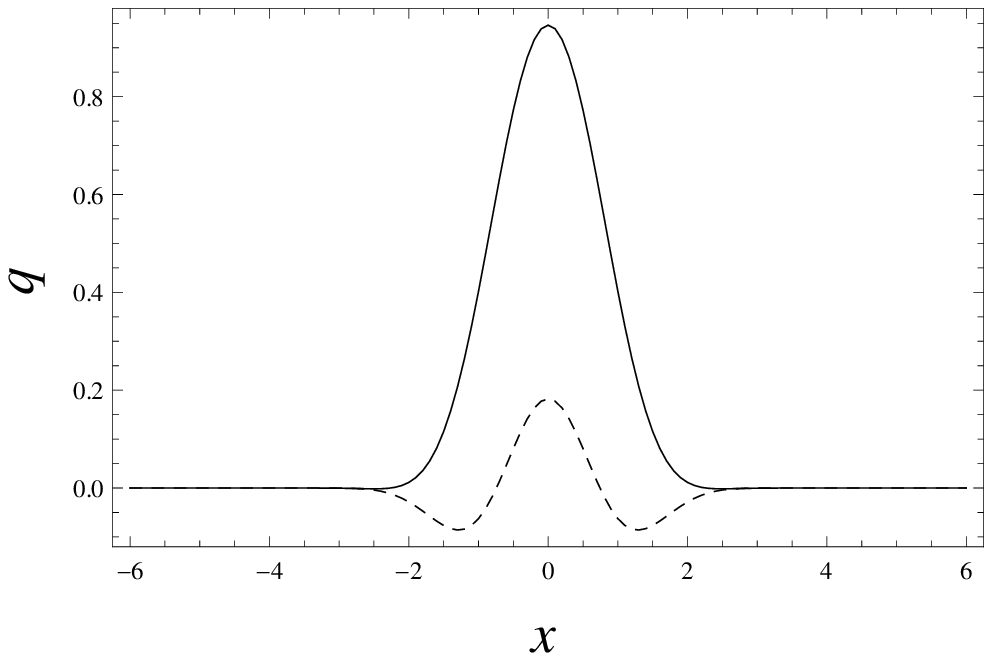}}
\subfigure[]{\includegraphics[width=.45\linewidth]{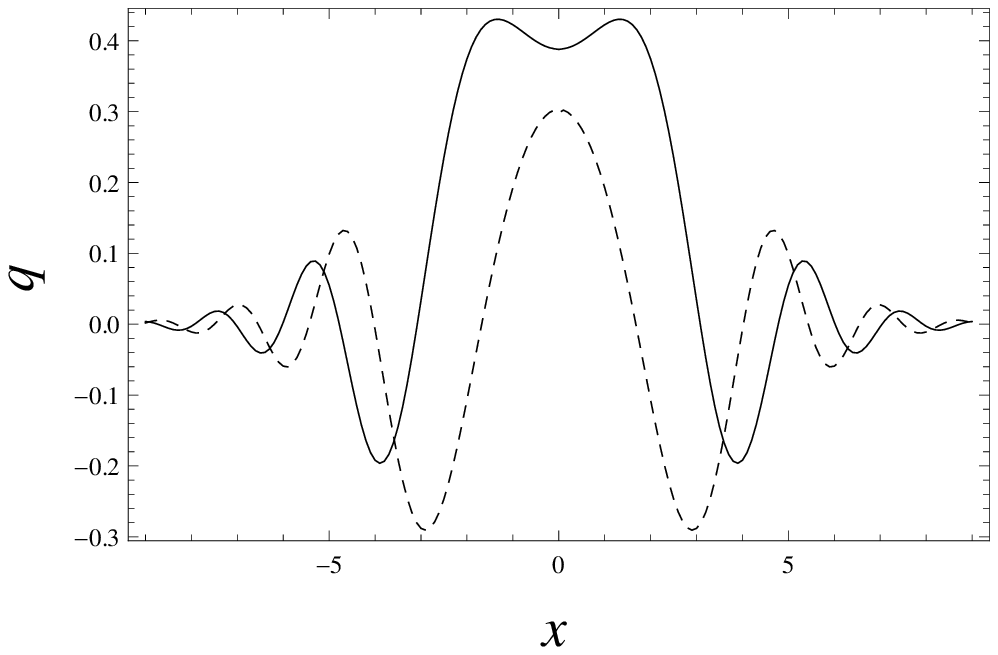}}
\caption{\label{FT-soln}(a)  A plot of the solution of \eqref{LS} at $t = 0.1$ (real part: solid, imaginary part: dashed).  (b)  A plot of the solution of \eqref{LS} at $t = 1$ (real part: solid, imaginary part: dashed).}
\end{figure}

\begin{remark}
When the initial condition $q_0$ has exponential decay the method of steepest descent for integrals combined with quadrature methods is the preferred way to solve this PDE \cite{TrogdonThesis}.
\end{remark}

\end{example}

\section{Oscillatory singular integral equations on the line}\label{Section:OscLine}

We have constructed an oscillatory basis of $L^2(\mathbb R)$ that is closed under function multiplication.  Additionally, the Cauchy operators leave the basis invariant.  Therefore if
\begin{align}\label{gen-form}
G(z)- I = \begin{mat} f_1(z) e^{i\alpha_1 z} & f_2(z)e^{i\alpha_2 z} \\ f_3(z)e^{i\alpha_3 z} & f_4(z)e^{i\alpha_4 z} \end{mat},
\end{align}
where each $f_i$ is expanded in the basis $\{R_{j}\}_{j=-N}^N$  then we may apply $\mathcal C[G;\mathbb R]$ to $G-I$.  We can also treat the case where $G-I$ is a sum of matrices of this form but this is beyond the scope of this paper.  Furthermore, we have a formula for the inner product.  These are all the required pieces to apply GMRES to $\mathcal C[G;\mathbb R]u = G-I$.  We discuss this in the examples below and we always use $\beta = 1$.

\subsection{Inverse scattering for small time}

We consider the numerical solution of \eqref{IST-RHP} where $\rho$ is the reflection coefficient associated with $q_0(x) = e^{-x^2}$.  We expand $\rho(z)e^{4ik^2t}$ in the basis $\{R_j\}_{j=-\infty}^\infty$ when $t$ is small.    The matrix $G(z;x,t)$ is of the form \eqref{gen-form} so that we may apply GMRES to \eqref{sie}.  As before only a fraction of the iterations of GMRES that are needed to solve $\mathcal C[G;\mathbb R]u=G-I$ are needed to solve $\mathcal C[G^{-1};\mathbb R]\mathcal C[G;\mathbb R]u = \mathcal C[G^{-1};\mathbb R](G-I)$.

When $t=0$ we are able to solve the preconditioned equation for moderate values of $x$.  There are two factors that must be taken into account when discussing efficiency of the method:
\begin{enumerate}
\item the number of GMRES iterations needed to reach a prescribed tolerance, and
\item the number of basis functions required to resolve the solution.
\end{enumerate}

\begin{figure}[tbp]
\centering
\includegraphics[width=.7\linewidth]{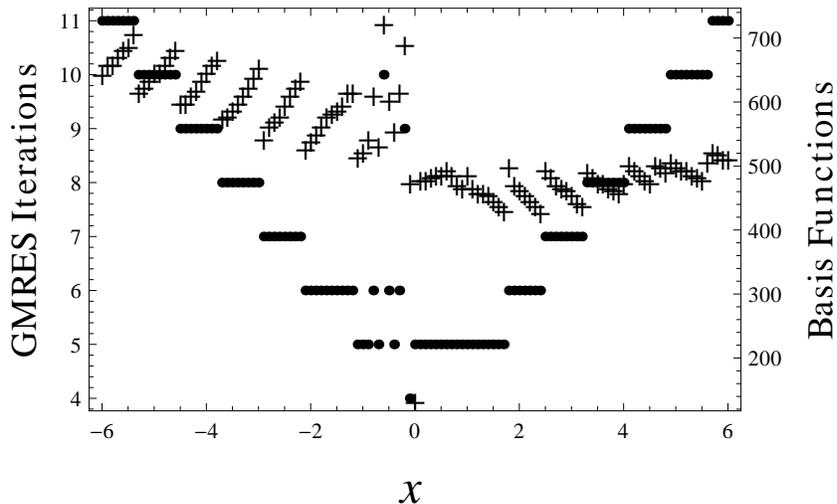}
\caption{\label{OscMatrixGMRES-reg}  A plot of the number of GMRES iterations (dots, left scale) to converge to a tolerance of $10^{-8}$ and number of needed basis functions (crosses, right scale) versus $x$.  We see that because we use oscillatory basis functions the number of needed basis functions appears to be bounded.  The operator is increasingly ill-conditioned as $|x|$ increases and we need more GMRES iterations.}
\end{figure}

While the number of basis functions to resolve the approximate solution at each stage of GMRES appears to be bounded, for moderate values of $x$ it takes many iterations of GMRES.  We display this behavior in Figure~\ref{OscMatrixGMRES-reg}.  In the following sections we present straightforward methods motivated by the Deift-Zhou method of nonlinear steepest descent \cite{DeiftZhouAMS} to reduce the number of GMRES iterations that are needed.

\subsection{Preconditioning for $x> 0$}\label{Section:precond-pos}

The matrix $G(z;x,t)$ in \eqref{IST-RHP} admits two important factorizations.  The first of which is
\begin{align}\label{MP}
\begin{split}
G(z;x,t) &= M(z;x,t)P(z;x,t),\\
P(z;x,t) &= \begin{mat} 1 & 0 \\ \rho(z)e^{2izx+4iz^2t} \end{mat},\\
M(z;x,t) &= \begin{mat} 1 & -\bar \rho(z) e^{-2izx-4iz^2t} \\ 0 & 1 \end{mat}.
\end{split}
\end{align}
The second factorization is discussed in the following section.  We rewrite \eqref{sie} using this factorization.  Write $G = MP$ and
\begin{align}
u - \mathcal C_{\mathbb R}^-u \cdot (MP-I) &= MP-I,\notag\\
uP^{-1} - \mathcal C_{\mathbb R}^- u\cdot(M-P^{-1}) &= M - P^{-1}. \label{precond-pos}
\end{align}

Before we discuss applying GMRES to \eqref{precond-pos}, we discuss the motivation for using such a factorization.  In the asymptotic analysis of solutions of the NLS equation with the Deift-Zhou method of nonlinear steepest descent the factorization \eqref{MP} is used for $t=0$ and $x\gg 0$.  The RHP \eqref{IST-RHP} on the line is deformed to one posed on two bi-infinite contours, parallel with real axis, with one lying in each of the upper- and lower-half planes.  The jump matrix on the upper contour is $P$ with $M$ being the jump on the lower contour.  A review of this can be found in \cite{deift-zhou:old-nls} (see also \cite{TrogdonSONLS}).  While in the present context we require no deformation off the real axis, we use this analysis to guide the choice \eqref{precond-pos}.

\begin{figure}[tbp]
\centering
\subfigure[]{\includegraphics[width=.48\linewidth]{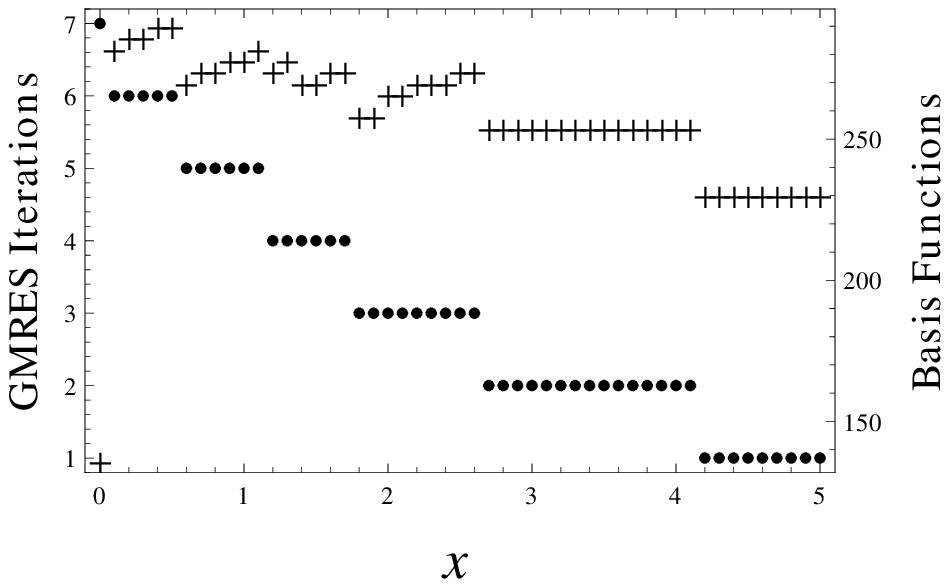}\label{OscMatrixGMRES-pos}}
\subfigure[]{\includegraphics[width=.48\linewidth]{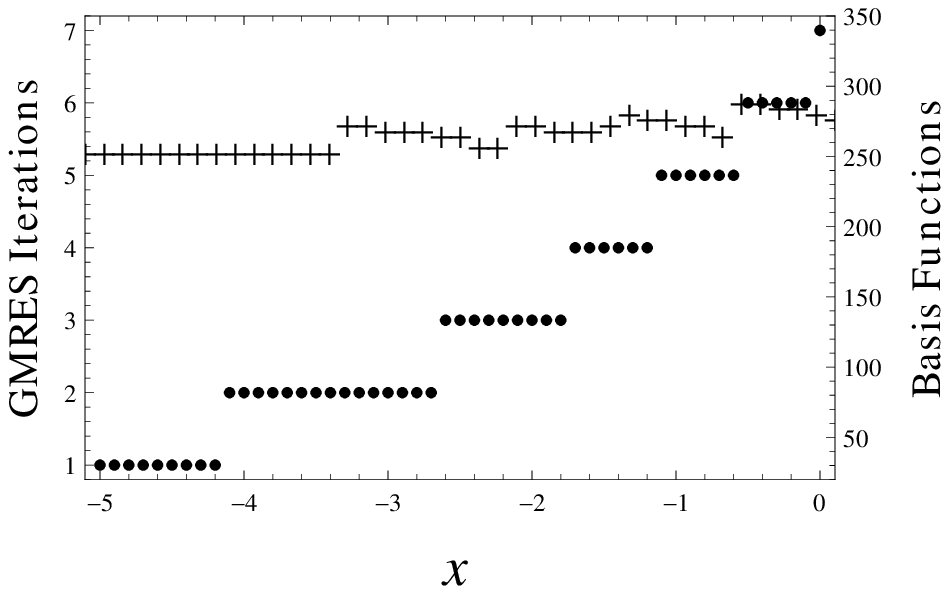}\label{OscMatrixGMRES-neg}}
\caption{A plot of the number of GMRES iterations (dots, left scale) to converge to a tolerance of $10^{-8}$ and number of needed basis function (crosses, right scale) versus $x$.  (a)  GMRES applied to \eqref{precond-pos}.  We see that because we use oscillatory basis functions the number of needed basis functions appears to be bounded and the number of GMRES iterations needed decreases as $x$ increases.  (b) GMRES applied to \eqref{precond-neg}.  The number of GMRES iterations needed decreases as $-x$ increases.}
\end{figure}

We apply GMRES directly to \eqref{precond-pos} without composition with any operator as in the preconditioned equation.  In practice, it takes so few iterations of GMRES to solve \eqref{precond-pos} that no further preconditioning is needed.  We demonstrate the efficiency of the computation in Figure~\ref{OscMatrixGMRES-pos}.  We find that fewer iterations of GMRES are needed for larger values of $x$ and the number of basis functions required is bounded.  This indicates bounded computational cost for all $x > 0$.

\subsection{Preconditioning for $x < 0$}\label{Section:precond-neg}

For $x< 0$ we take a similar approach as in \eqref{Section:precond-pos} and factor the matrix $G$.  In this case the matrix factorization is more complicated.  First, we note that
\begin{align*}
G(z;x,t) &= L(z;x,t)D(z)U(z;x,t),\\
L(z;x,t) &= \begin{mat} 1 & 0 \\
\displaystyle\frac{\rho(z)e^{2ixz+4iz^2t}}{1-|\rho(z)|^2} & 1 \end{mat},\\
D(z) &= \begin{mat} 1-|\rho(z)|^2 & 0 \\ 0 & (1-|\rho(z)|^2)^{-1} \end{mat},\\
U(z;x,t) &= \begin{mat} 1 & \displaystyle-\frac{\bar \rho(z) e^{-2ixz+4iz^2t}}{1-|\rho(z)|^2} \\ 0 & 1 \end{mat}.
\end{align*}
The matrix $D$ admits a Riemann--Hilbert factorization:
\begin{align*}
\Delta^+(z) = \Delta^-(z) D(z),~~ \Delta(z) = \diag (\delta(z), \delta^{-1}(z)),\\
\delta(z) = \exp\left(\frac{1}{2 \pi i} \int_{\mathbb R} \frac{ \log(1-|\rho(s)|^2)}{s-z} ds   \right).
\end{align*}

Next, we note that since \eqref{integral-reconstruct} involves only the $(2,1)$-component of $u$ (and of $\Phi$) we may consider
\begin{align*}
\tilde G(z;x,t) &= \Delta^-(z)L(z;x,t)D(z)U(z;x,t)(\Delta^+)^{-1}(z)\\
&= \Delta^-(z)L(z;x,t)(\Delta^-)^{-1}(z) \Delta^+(z) U(z;x,t) (\Delta^+)^{-1}(z)\\
& = \tilde L(z;x,t) \tilde U(z;x,t),
\end{align*}
where
\begin{align*}
\tilde U(z;x,t) = \begin{mat} 1 & \displaystyle -\frac{\bar \rho(z) e^{-2ixz+4iz^2t}}{1-|\rho(z)|^2}(\delta^+)^2(z) \\ 0 & 1 \end{mat}, ~~ \tilde L(z;x,t) = \begin{mat} 1 & 0 \\
\displaystyle\frac{\rho(z)e^{2ixz+4iz^2t}}{1-|\rho(z)|^2}(\delta^-)^{-2}(z) & 1 \end{mat}.\\
\end{align*}

We apply GMRES to the operator equation
\begin{align}
u\tilde U^{-1} - \mathcal C_{\mathbb R}^- u\cdot(\tilde L-\tilde U^{-1}) &= \tilde L - \tilde U^{-1}. \label{precond-neg}
\end{align}
We do not consider preconditioning it further. The functions $\delta^\pm(z)$ may be computed with the techniques in Section~\ref{Section:scalar}.   We demonstrate the accuracy and efficiency of the method in Figure~\ref{OscMatrixGMRES-neg}.  As in the case  of positive $x$ we see that the number of GMRES iterations needed to converge decays as $-x$ increases.  Furthermore, the number of basis functions needed to accurately resolve the solution appears to be bounded.  This indicates bounded computational cost for all $x$.

\begin{figure}[tbp]
\centering
\subfigure[]{\includegraphics[width=.49\linewidth]{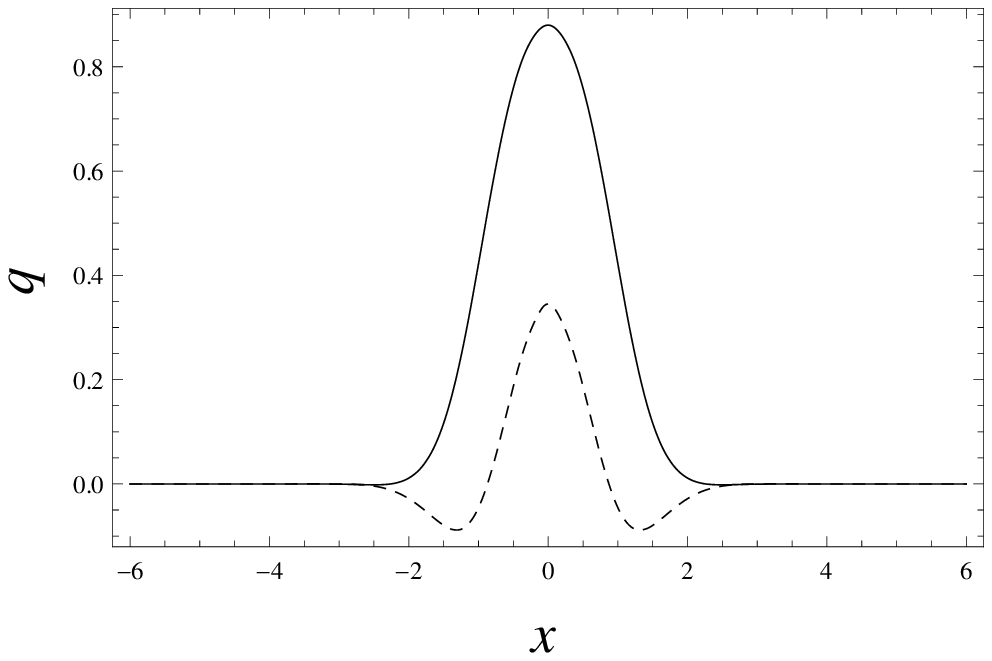}}
\subfigure[]{\includegraphics[width=.49\linewidth]{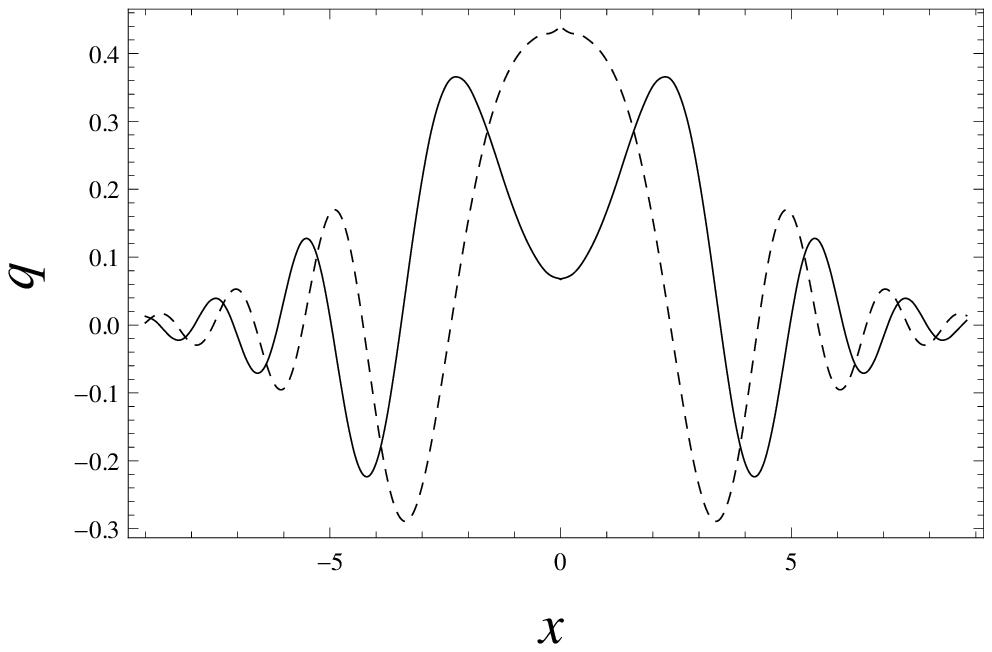}}
\caption{\label{NLS-1} (a) A plot of the solution of \eqref{NLS} at $t = 0.1$ (real part: solid, imaginary part: dashed). (b)  A plot of the solution of \eqref{NLS} at $t = 1$ (real part: solid, imaginary part: dashed).}
\end{figure}

We combine the approach for $x>0$ in Section~\ref{Section:precond-pos} with the approach for $x<0$ in this section to compute the solution of the NLS equation for $t > 0$. The approach is limited to small $t$ and we show the solution of the NLS equation with $q_0(x)=e^{-x^2}$ when $t = 0.1, 1.0$ in Figure~\ref{NLS-1}.

\subsection{Jump matrices with slow decay}

\begin{figure}[tbp]
\centering
\includegraphics[width=.5\linewidth]{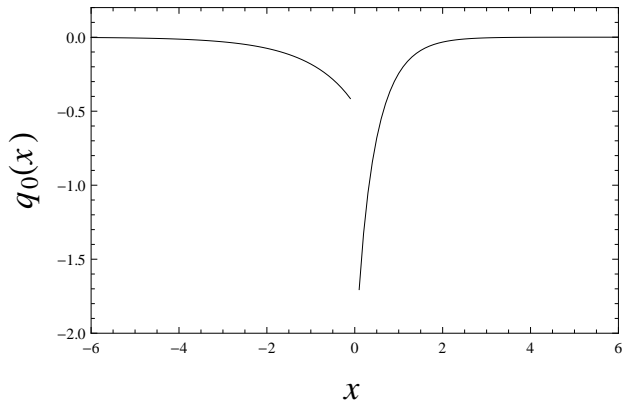}
\caption{\label{discontinuous} A plot of a discontinuous function computed by solving \eqref{sie} and evaluating \eqref{integral-reconstruct} when $\rho(z) = .9i/(z-i)$.}
\end{figure}

Slow decay in the reflection coefficient indicates a lack of smoothness (even lack of continuity!) in the initial condition.  The method presented here is well-suited to deal with slow decay.  We present an example of this.  Let $\rho(z) = 0.9i/(z-i)$.  This function is easily represented in the basis $\{R_j\}$:
\begin{align*}
\rho(z) = 0.45 R_{-1}(z).
\end{align*}
As before, we combine the approach for $x>0$ in Section~\ref{Section:precond-pos} with the approach for $x<0$ in Section~\ref{Section:precond-neg} to solve \eqref{sie} and evaluate \eqref{integral-reconstruct} at $t=0$.  See Figure~\ref{discontinuous} for a plot of the initial condition.

\begin{remark}
Due to the slow decay of $\rho(z)$, $\rho(z)e^{4iz^2t}$ cannot efficiently be represented by the basis $\{R_j\}$.  Indeed, the first derivative of $\rho(z)e^{4iz^2t}$ with respect to $z$ does not decay at infinity and the fast Fourier transform based technique presented above fails.
\end{remark}

\section{Conclusions}

We have constructed a new numerical method for the solution of oscillatory singular integral equations.  We derived formulae for the action of the Cauchy operators, integration and inner products and function multiplication.  This allows us to apply the infinite-dimensional GMRES algorithm to singular integral equations on the real axis.  In the examples, we explored computing the inverse scattering transform for small time.  While the method in it current state does not beat the state of the art \cite{TrogdonSONLS,TrogdonSOKdV} in terms of speed and scope it does have some important implications for future research:
\begin{itemize}
\item the Fredholm regulator of a singular integral operator provides an effective GMRES preconditioner,
\item matrix factorizations used in the asymptotic analysis of RHPs also provide effective preconditioners,
\item the method allows for slow decay in the reflection coefficient $\rho$ which allows us to perform inverse scattering for discontinuous potentials at $t = 0$, and
\item the effectiveness of GMRES in this context gives a strong indication that the method could be expanded to allow for the effective solution of oscillatory singular integral equations with more complicated oscillatory basis functions, \emph{i.e.},  bases that include the factors $e^{4iz^2t}$ and $e^{8iz^3t}$.   This may result the numerical solution of RHPs that arise in the inverse scattering transform for all $x$ and $t$ without deformation.
\end{itemize}

Specifically, for a scalar singular integral equation,  using the Fredholm regulator as a preconditioner reduced the number of GMRES iterations from 20 to 4 to achieve the same tolerance (see Figure \ref{ScalarGMRES-reg}).  For a matrix singular integral equation, preconditioning with the Fredholm regulator reduced the number of GMRES iterations from 34 to 6 to achieve the same tolerance (see Figure~\ref{MatrixGMRES-reg}). When considering oscillatory singular integral equations we use simple algebraic preconditioners that can reduce the required number of GMRES iterations by more than 90\% (see Figures~\ref{OscMatrixGMRES-neg} and \ref{OscMatrixGMRES-pos}).

\section*{Acknowledgments}

We acknowledge the National Science Foundation for its generous support through grants NSF-DMS-1008001 and NSF-DMS-1303018.  Any opinions, findings, and conclusions or recommendations expressed in this material are those of the authors and do not necessarily reflect the views of the funding sources. We also thank the anonymous referees for their input which improved this manuscript.

\bibliographystyle{plain}
\bibliography{GMRESstyle}

\end{document}